\documentclass[11pt]{article}

\usepackage[a4paper,margin=1in]{geometry}
\usepackage{amsmath,amssymb,amsthm,mathtools}
\usepackage{graphicx}
\usepackage{float}
\usepackage[expansion=false]{microtype}
\usepackage{etoolbox}
\usepackage[hidelinks]{hyperref}
\ifdefined\pdfinfoomitdate\pdfinfoomitdate=1\fi
\ifdefined\pdftrailerid\pdftrailerid{}\fi
\ifdefined\pdfsuppressptexinfo\pdfsuppressptexinfo=15\fi

\apptocmd{\normalsize}{%
  \setlength{\abovedisplayskip}{5pt plus 2pt minus 3pt}%
  \setlength{\belowdisplayskip}{5pt plus 2pt minus 3pt}%
  \setlength{\abovedisplayshortskip}{2pt plus 1pt minus 1pt}%
  \setlength{\belowdisplayshortskip}{2pt plus 1pt minus 1pt}}{}{}
\apptocmd{\small}{%
  \setlength{\abovedisplayskip}{4pt plus 2pt minus 2pt}%
  \setlength{\belowdisplayskip}{4pt plus 2pt minus 2pt}%
  \setlength{\abovedisplayshortskip}{2pt plus 1pt minus 1pt}%
  \setlength{\belowdisplayshortskip}{2pt plus 1pt minus 1pt}}{}{}
\normalsize
\apptocmd{\thebibliography}{\setlength{\itemsep}{1pt plus .5pt}}{}{}

\makeatletter
\pretocmd{\@citex}{%
  \@for\paper@citekey:=#2\do{%
    \edef\paper@citekey{\expandafter\@firstofone\paper@citekey\@empty}%
    \expandafter\gdef\csname paper@cited@\paper@citekey\endcsname{}%
  }%
}{}{\PackageError{paper}{Could not enable citation filtering}{}}
\newcommand{\citedbibitem}[2]{%
  \ifcsname paper@cited@#1\endcsname
    \bibitem{#1}#2%
  \fi
}
\makeatother

\newcommand{\T}{{\mathbb{T}^{3}}}
\newcommand{\R}{{\mathbb{R}^{3}}}
\newcommand{\Z}{{\mathbb{Z}^{3}}}
\newcommand{\dd}{\,\mathrm d}
\newcommand{\norm}[1]{\left\lVert #1\right\rVert}
\newcommand{\one}{\mathbf 1}
\theoremstyle{plain}
\newtheorem{theorem}{Theorem}[section]
\newtheorem{proposition}[theorem]{Proposition}
\newtheorem{lemma}[theorem]{Lemma}
\newtheorem{corollary}[theorem]{Corollary}
\theoremstyle{definition}
\newtheorem{definition}[theorem]{Definition}
\theoremstyle{remark}
\newtheorem{remark}[theorem]{Remark}

\title{Hypocoercivity
for a thermostatted Boltzmann equation}
\author{Josephine Evans \and Alexander Pilakoutas}
\date{}

\hypersetup{pdftitle={Hypocoercivity for a thermostatted Boltzmann equation},
  pdfauthor={Josephine Evans, Alexander Pilakoutas}}
\begin{document}
\maketitle
\begin{abstract}
We study the Boltzmann equation for a gas on the three-dimensional
torus. Its particles collide with one another and with a thermal reservoir
whose temperature varies in space.
The reservoir heats the gas in some regions and cools it in others, driving a
nonequilibrium energy current.

For the isotropic cutoff kernel, with Maxwell molecules or hard potentials,
we prove global
well-posedness and exponential convergence to the same stationary
state for initial data with a common Gaussian upper bound and
Gaussian moment bound. These bounds determine how weak self
collisions must be relative to bath collisions. The estimates
are uniform over bath profiles between fixed positive temperature
bounds.

For any nonconstant reservoir profile on the three-dimensional torus,
this stationary state cannot be a local Maxwellian.  When particles collide only with
the bath, relative entropy to the stationary state cannot increase.  Once
self collisions are added, arbitrarily small
perturbations can initially increase it, even though they later return
exponentially.

\end{abstract}


\section{Introduction}
This paper is about nonequilibrium steady states (NESS) in kinetic theory.
We wish to understand steady states in gases which interact with thermal
reservoirs, meaning the steady state can maintain a constant flux of
energy. We will study a gas modelled by the Boltzmann equation interacting
with thermostat terms. When the reservoir temperature is nonconstant,
the steady states will not be Maxwellians. In fact, a priori, the question
of existence and uniqueness of these states is open. The question of
stability of these steady states is mathematically challenging. They are
non-explicit and the traditional tools are designed for working around
Maxwellians. For sufficiently weak self collisions, we prove global well-posedness and exponential convergence to a unique stationary state in a prescribed class with Gaussian pointwise and moment bounds. The initial data need not be close to the stationary state; the permitted strength of self collisions depends on these bounds. The results include hard-sphere collisions and allow arbitrary reservoir temperature variation between fixed positive bounds. For a nonconstant reservoir temperature, the stationary state is not a local Maxwellian. Despite exponential relaxation, arbitrarily small perturbations can initially increase relative entropy to this state. A key motivation is understanding the emergence of
heat conduction in gases. Once we have a steady state, we can see a
relationship between its temperature and heat flux which should converge to
Fourier's law in an appropriate scaling limit.

\subsection{Model and Notation}

We study
\begin{equation}\label{eq:boltzmann}
  \partial_t f+v\cdot\nabla_x f
  =(1-\alpha)Q(f,M_{\tau(x)})+\alpha Q(f,f),
  \qquad (x,v)\in\mathbb T^3\times\mathbb R^3.
\end{equation}
Here \(f(t,x,v)\) is the gas density and \(0\leq\alpha<1\).
We work with probability densities on phase space, including
initial data and stationary states. The terms \(Q(f,f)\) and
\(Q(f,M_{\tau(x)})\) describe self collisions and collisions with
a Maxwellian reservoir at temperature \(\tau(x)\), respectively. In \(Q\), collisions are weighted by \(|v-v_*|^\gamma\), \(\gamma\geq0\);
see Section~\ref{sec:models} where we describe the models and in particular $Q$ in more detail.
We assume
\[
  0<\tau_-\leq\tau(x)\leq\tau_+<\infty
  \qquad\text{for almost every }x\in\mathbb T^3.
\]

For probability densities \(f\) and \(F>0\), the relative
entropy of \(f\) with respect to \(F\) is
\[
  H(f\,|\,F)=\int_{\mathbb T^3\times\mathbb R^3}
  \left(f\log\frac fF-f+F\right)\,\dd v\dd x .
\]
We also use the weighted norms
\[
 \norm{h}_{L^1(W)}=\int_{\T}\int_{\R}|h(x,v)|W(v)\,\dd v\dd x,
 \qquad
 \norm{h}_{L^\infty(W)}=
 \operatorname*{ess\,sup}_{(x,v)\in\T\times\R}|h(x,v)|W(v).
\]
We write \(\langle v\rangle=(1+|v|^2)^{1/2}\) and
\(W(v)=e^{a|v|^2}\), with \(0<a<1/(2\tau_+)\).

\subsection{Main Results}

When collisions with the reservoir dominate, we show the existence of a NESS as well as exponential convergence towards it. When \(\tau\) is not constant and \(\alpha >0\), then the NESS is not a local Maxwellian,
and relative entropy can increase.

\begin{theorem}[Relaxation to the stationary state]
  \label{thm:intro-relaxation}\label{thm:boltzmann-contraction}
  There are \(R_*,t_c>0\), depending only on
  \(a,\gamma,\tau_-,\tau_+\), such that, for every \(R\geq R_*\),
  there are \(C>0\) and \(\alpha_0\in(0,1)\), depending only on
  these parameters and \(R\), for which the following holds.
  The class
  \[
   \mathcal R=\left\{f\geq0:\quad
   \int_{\T}\int_{\R}f\,\dd v\dd x=1,\quad
   \|f\|_{L^1(W)}+\|f\|_{L^\infty(W)}\leq R\right\}
  \]
is nonempty. For \(0\leq\alpha\leq\alpha_0\) and initial datum
\(f_{\rm in}\in\mathcal R\), equation~\eqref{eq:boltzmann} has a unique
global mild solution \(f\). The solution map at time \(t_c\) is a strict
contraction from \(\mathcal R\) into itself for the sum of the two weighted
norms. There is a unique stationary state \(F_\alpha\) in \(\mathcal R\),
and the solution satisfies
  \begin{align*}
   &\|f(t)-F_\alpha\|_{L^1(W)}+\|f(t)-F_\alpha\|_{L^\infty(W)}\\
   &\quad+\int_t^\infty\|\langle v\rangle^\gamma(f(s)-F_\alpha)\|_{L^1(W)}\,\dd s\\
   &\qquad\leq C\,2^{-\lfloor t/t_c\rfloor}
   \bigl(\|f_{\rm in}-F_\alpha\|_{L^1(W)}
   +\|f_{\rm in}-F_\alpha\|_{L^\infty(W)}\bigr),\qquad t\geq0.
  \end{align*}
\end{theorem}

The upper bound \(1/(2\tau_+)\) for the Gaussian exponent \(a\)
is sharp when \(\tau_+=\operatorname{ess\,sup}\tau\).
The stationary state cannot satisfy \(F_\alpha\le Ce^{-a|v|^2}\)
with \(a>1/(2\tau_+)\); see Remark~\ref{rem:stationary-tails}.

While we don't expect our theorem to give a sharp condition on $\alpha$ for convergence to a steady state, some constraints on \(\alpha\) may be needed for such a result to hold. In some cases, temperature gradients may drive movement.
For instance, in the presence of gravity, a gas that is heated from below may experience convection;
this has been analyzed close to the hydrodynamic limit of the Boltzmann equation \cite{ArkerydEspositoMarraNouri2010,ArkerydEspositoMarraNouri2012}.
Gravity is not present in \eqref{eq:boltzmann}, but temperature gradients can amplify sound waves.
As Lord Rayleigh explained \cite[\S~322g, p.~226]{Rayleigh1896}, ``If heat be given to
the air at the moment of greatest condensation, or be taken from it at the
moment of greatest rarefaction, the vibration is encouraged.''
The same mechanism drives the Sondhauss and Rijke tubes
\cite[Sections~322i--322j, pp.~230--234]{Rayleigh1896}.
For a review of thermoacoustic engines, see \cite{Swift1988}.

When \(\alpha=0\), the evolution operator forms a Markov semigroup, so the relative entropy cannot increase
\cite{Csiszar1967,Morimoto1963}. The next theorem shows that self
collisions can cause an initial increase, even arbitrarily close to the stationary state.

\begin{theorem}[Entropy growth]
\label{thm:intro-entropy-growth}
In addition to the hypotheses of
Theorem~\ref{thm:intro-relaxation}, suppose that
\(\tau\) is not almost everywhere constant and
\(0<\alpha<\alpha_0\).  Then \(Q(F_\alpha,F_\alpha)\not\equiv0\).
For every \(\delta>0\), there is a datum
\(f_{\rm in}\in\mathcal R\) such that
\[
 \|f_{\rm in}-F_\alpha\|_{L^1(W)}
 +\|f_{\rm in}-F_\alpha\|_{L^\infty(W)}<\delta,
 \qquad \frac12\leq\frac{f_{\rm in}}{F_\alpha}\leq\frac32,
\]
and the global solution satisfies
\[
 H(f(t)\,|\,F_\alpha)>H(f_{\rm in}\,|\,F_\alpha)
 \qquad\text{for all sufficiently small }t>0.
\]
The same solution satisfies the exponential return estimate in
Theorem~\ref{thm:intro-relaxation}.
\end{theorem}

The calculation in Proposition~\ref{prop:entropy-three-terms} uses the
same changes of variables as the H-theorem. For smooth solutions, it
decomposes the derivative of relative entropy to \(F\) into dissipation
from reservoir collisions,
dissipation from self collisions, and a cross term from self collisions.
The dissipations cannot increase relative entropy. The cross term can
have either sign and vanishes when \(F\) is a local Maxwellian, which
occurs only when \(\tau\) is constant
(Lemma~\ref{lem:stationary-maxwellian-rigidity}). For nonconstant \(\tau\)
and \(\alpha>0\), we construct arbitrarily small perturbations of \(F\)
that change both the density and velocity distribution, giving solutions
whose relative entropy exceeds its initial value for small positive times
(Section~\ref{sec:entropy-growth}).

\subsection{Previous Work}
\label{sec:intro-comparison}

\subsubsection{Nonequilibrium Stationary States}
\label{sec:intro-states}

A useful survey of nonequilibrium steady states and their hydrodynamic limits in various kinetic models is Esposito and Marra
\cite{EspositoMarra2020}.

\paragraph{Boltzmann models.}
For stationary heat conduction, many results concern a small variation
in the wall temperature or the limit of small Knudsen number
\cite{EspositoGuoKimMarra2013,EspositoGuoKimMarra2018}.

For diffuse reflection at wall temperature \(\theta\), Esposito, Guo,
Kim and Marra \cite{EspositoGuoKimMarra2013} constructed a locally unique,
dynamically stable stationary state for small \(|\theta-\theta_0|\).
Arkeryd and Nouri \cite{ArkerydNouri2000} obtained weak \(L^1\) stationary
solutions in a slab with diffuse reflection, without requiring small
variation of \(\theta\). Arkeryd \cite{Arkeryd2000} treated a convex domain
with particles entering through the boundary and a collision operator
truncated at small velocities.

Stochastic coupling to reservoirs goes back to Bergmann and Lebowitz
\cite{BergmannLebowitz1955}.
For \(\alpha=0\) and constant temperature, \eqref{eq:boltzmann} is the linear Boltzmann
equation, for which exponential convergence to equilibrium is known in
the spatially homogeneous case \cite{BisiCanizoLods2015} and on the torus
\cite{CanizoCaoEvansYoldas2020}. For a spatially uniform thermal medium,
Fr\"ohlich and Gang \cite{FroehlichGang2012} studied a nonlinear
Boltzmann equation on \(\T\), with a reservoir transition kernel
satisfying detailed balance. Self collisions carry a small
coefficient playing the role of \(\alpha\). Under their kernel
assumptions, they proved exponential convergence in \(L^1\) to a Maxwellian
for sufficiently smooth data with polynomial velocity moments. Here the reservoir
temperature can vary in space.

Even in spatially homogeneous models, reservoirs at different temperatures
can give a stationary distribution that is not Maxwellian. For spatially homogeneous cutoff Maxwell molecules, Carlen, Lebowitz
and Mouhot \cite{CarlenLebowitzMouhot2015} proved uniqueness and
exponential convergence in a Fourier probability metric.
For Kac's model with a non-Maxwellian thermostat, Evans \cite{Evans2016}
proved convergence to a unique stationary state, with estimates as the
particle number tends to infinity.
For granular gases, rescaling around a self-similar profile leads to a
stationary problem. Mischler and Mouhot proved existence of these profiles
\cite{MischlerMouhot2006} and their uniqueness and stability for small
inelasticity \cite{MischlerMouhot2009}.

\paragraph{BGK models.}
In the BGK equation, binary self collisions are replaced by
relaxation towards a Maxwellian with the local density, mean
velocity and temperature of the gas.
For the spatially homogeneous thermostats studied by Carlen and
collaborators, the steady distribution can be written explicitly.
Their work on kinetic Fokker--Planck and BGK equations on the torus
\cite{CarlenEspositoLebowitzMarraMouhot2018} gives exponential relaxation
for spatially uniform data and uniqueness for some simplified models.
For the one-dimensional BGK equation with reservoirs at different
temperatures, the small-\(\alpha\) uniqueness theorem assumes finite
energy and spatial density in \(L^p\), \(p>1\). Local stability holds
throughout \(0\le\alpha<1\)
\cite{CarlenEspositoLebowitzMarraMouhot2020}.

In \cite{EvansMenegaki2021}, Evans and Menegaki constructed a BGK
steady state between thermal walls by Schauder's theorem, assuming large
and well-separated wall temperatures. Their later torus model
\cite{EvansMenegaki2025} is the BGK counterpart of
\eqref{eq:boltzmann}. For a sufficiently large lower bound on the
reservoir temperature and small \(\alpha\), they obtained uniqueness
among steady states satisfying temperature bounds and a spectral
gap for the linearisation. Nonlinear stability was left open
\cite[Section~1.3.2]{EvansMenegaki2025}. We prove nonlinear stability
for the Boltzmann model under the bounds \(0<\tau_-\le\tau\le\tau_+<\infty\).

\paragraph{Other kinetic models.}
For free transport with Maxwell boundary conditions, convergence to a
steady state is possible without a small temperature difference
\cite{Bernou2025}. For the kinetic Fokker--Planck equation with
constant wall temperature, Carrapatoso and Mischler
\cite{CarrapatosoMischler2024} obtain weighted \(L^1\)-to-\(L^\infty\)
bounds and exponential convergence. Carrapatoso, Gabriel, Medina and
Mischler \cite{CarrapatosoGabrielMedinaMischler2024} allow spatially
varying wall temperatures; in the mass-conserving case, their result
gives a unique steady probability density and exponential convergence. A weakly nonlinear
Fokker--Planck equation with BGK thermostats is studied in
\cite{EvansMedina2026}, where stability holds for small perturbations
of a steady state.

\subsubsection{Well-Posedness and Stability}
\label{sec:intro-stability}

\paragraph{Well-posedness.}
Global renormalised solutions are given by the theorem of DiPerna and Lions \cite{DiPernaLions1989}. Here we use bounded
mild solutions, whose local theory with Gaussian velocity bounds goes
back to Kaniel and Shinbrot \cite{KanielShinbrot1978}.
Global results for bounded solutions include the near-vacuum problem
in the whole space \cite{IllnerShinbrot1984} and perturbations of a
Maxwellian \cite{GualdaniMischlerMouhot2017,Guo2010,Ukai1974}.
For background on classical Boltzmann solution theory, see
\cite{CercignaniIllnerPulvirenti1994,Glassey1996,Villani2002Handbook}.
We allow data far from the
stationary state, with the allowed size of \(\alpha\) depending on
their bound in \(L^1(W)\cap L^\infty(W)\).

\paragraph{Spectral and hypocoercive methods.}
Near a global Maxwellian \(M\), stability proofs combine the spectral
gap of the linearised collision operator in \(L^2(M^{-1})\) with estimates
for transport \cite{GualdaniMischlerMouhot2017,Guo2010,Ukai1974}.
For a stationary state \(F\), the quadratic term in \(H(F+h\mid F)\)
is \(\frac12\int|h|^2/F\,\dd v\dd x\), suggesting the space
\(L^2(F^{-1})\). The linearised flow need not decrease this norm.
It does for the uniform BGK stationary state in
\cite{CarlenEspositoLebowitzMarraMouhot2020}. Evans and Menegaki
\cite[Section~1.3.2]{EvansMenegaki2025} discuss growth of this norm for
the thermal-wall problem. For their torus model, they proved exponential
decay using a modified quadratic energy
\cite[Sections~1.3.1 and~3.6]{EvansMenegaki2025}, following Dolbeault,
Mouhot and Schmeiser \cite{DolbeaultMouhotSchmeiser2015}.
For nonconstant \(\tau\) and \(0<\alpha<\alpha_0\), the quadratic form
\eqref{eq:linearised-three-terms} in our entropy-growth proof has a
positive direction.

\paragraph{Probabilistic methods and weak nonlinear coupling.}
Harris's theorem estimates convergence from properties of the transition
kernel, without a formula for the stationary state. Ca\~nizo, Cao, Evans
and Yolda\c{s} \cite{CanizoCaoEvansYoldas2020} used it for the linear
Boltzmann equation on the torus, and Bernou \cite{Bernou2025} for linear
kinetic equations with thermal walls. A related argument using Doeblin's
theorem appears in \cite{CarlenEspositoLebowitzMarraMouhot2018}.
For spatially homogeneous equations, other contraction estimates use the
Gabetta--Toscani--Wennberg metric or Wasserstein distances
\cite{CarlenLebowitzMouhot2015,Evans2016}.

The nonlinear argument has precedents in the theory of
McKean--Vlasov diffusions \cite{Butkovsky2014} and jump processes whose
rates and kernels depend on their law \cite{Monmarche2023}. For a run-and-tumble equation, Evans and Yolda\c{s}
\cite{EvansYoldas2023} estimate the nonlinear term in Duhamel's formula
using the decay supplied by Harris's theorem. These results allow data
far from stationarity under moment bounds and a small coupling condition.
For our equation, the estimates must account for the dependence on
\(f(x,\cdot)\) and, for hard potentials, the growth of the collision
rate with velocity.

\paragraph{Transfer between norms.}
Gualdani, Mischler and Mouhot \cite{GualdaniMischlerMouhot2017}
developed a factorisation method for transferring semigroup decay to
larger spaces, including spaces with \(L^1\) integrability in velocity
for the Boltzmann equation. Mischler and Mouhot
\cite{MischlerMouhot2016} also treat transfer to smaller spaces.
One part of the generator gives a dissipative evolution; compositions
with the other improve integrability when Duhamel's formula is iterated. In the run-and-tumble equation,
Mischler and Weng \cite[Lemmas~3.4--3.6]{MischlerWeng2017} expand
Duhamel's formula and estimate the repeated-gain terms through
\(L^1\), bounding the remaining terms in \(L^\infty\). They obtain
bounds in a weighted \(L^1\cap L^\infty\) space. We use this
decomposition to extend Harris's \(L^1(W)\) decay to \(L^\infty(W)\).
The latter bound is needed to estimate the nonlinear collision term
(Lemma~\ref{lem:gain-bounds}).

\subsubsection{Relative Entropy}
\label{sec:intro-entropy}

For the Boltzmann equation on the torus without reservoirs, relative
entropy to global Maxwellian equilibrium is nonincreasing by the
H-theorem. In the inhomogeneous equation, collision dissipation vanishes
on local Maxwellians, so quantitative convergence also uses transport
\cite{DesvillettesVillani2005,Villani2009}. For the spatially homogeneous linear
Boltzmann equation, Bisi, Ca\~nizo and Lods
\cite{BisiCanizoLods2015} bounded relative entropy by a multiple of its
dissipation.

For nonlinear equations with a non-Maxwellian stationary state, the
classical H-theorem does not determine the sign of the derivative of
relative entropy to that state; see the discussion in
\cite{CarlenEspositoLebowitzMarraMouhot2020}.
A related failure of monotonicity occurs in the Vlasov--Fokker--Planck
model of Cesbron and Herda \cite{CesbronHerda2024}. Their usual free energy
can initially increase when the interaction potential is not even
\cite[Corollary~3.6]{CesbronHerda2024}.
Monmarch\'e \cite{Monmarche2025} constructs a different, nonincreasing
free energy from relative entropy per particle for the associated
interacting system.

Proposition~\ref{prop:stationary-entropy-balance}
accounts for the entropy of the gas and reservoir together and shows
strictly positive total entropy production at stationarity when
\(\tau\) is nonconstant. For nonconstant \(\tau\) and \(\alpha>0\),
we do not know how to correct \(H(f\,|\,F_\alpha)\) explicitly so that
it cannot increase along the nonlinear flow and is minimised only at
\(F_\alpha\).

\subsection{Method of Proof}
\label{sec:intro-method}

The key fact that makes our proofs work is good control on the linear semigroup (when $\alpha=0$). The linear semigroup is in a nonequilibrium setting, meaning that the usual entropy estimates around a Maxwellian do not directly apply. However, since it is a Markov process we are able to use Harris's theorem to get contraction in a weighted $L^1$ space. We then want to consider the equation with $\alpha>0$ but small as a perturbation of this equation. In order to control this perturbation we need weighted $L^\infty$ control on the linear semigroup. We get this using operator splitting arguments following ideas of Mischler \cite{MischlerWeng2017}.

To transfer the $L^1(W)$ decay to $L^\infty(W)$, we expand Duhamel's formula twice. This allows us to change variables from an intermediate velocity to an earlier collision position. The resulting term is bounded by an $L^1(W)$ contribution and a small multiple of the $L^\infty(W)$ bound. The terms containing the initial data decay in $L^\infty(W)$. Absorbing the small terms and using Harris's theorem gives decay in $L^\infty(W)$. For hard potentials, the collision rate grows with velocity. This introduces higher velocity moments, whose integrals in time we control.

\paragraph{Local nonlinear solutions.}
We first construct solutions to the nonlinear problem locally in time. We do this through a contraction mapping argument based on freezing the nonlinearity. The time interval $[0,T]$ on which this argument works depends on $\alpha$ and the initial bounds. We then want to show we can continue these solutions globally in time and that any two solutions contract together. If we begin with initial data in the bounded class $\mathcal{R}$ of Theorem~\ref{thm:intro-relaxation}, with sufficiently large radius, we first choose a time $t_c$ using the linear contraction. We then choose $\alpha$ small enough depending on $\mathcal{R}$ so that the solutions exist up to $t_c$ and return to $\mathcal{R}$ at that time. This allows us to patch together solutions to create a global solution.
For two solutions, their difference \(w=f-g\) satisfies
\begin{equation}\label{eq:nonlinear-difference-duhamel}
 w(t)=P_t^\tau w(0)+\alpha\int_0^t P_{t-s}^\tau
 \bigl[Q(w,f)+Q(g,w)-Q(w,M_{\tau(\cdot)})\bigr](s)\,\dd s.
\end{equation}
The weighted $L^1$ bounds on $Q$ use a Gaussian pointwise bound on one factor and a weighted moment of the other. Estimates for the linear equation with a source term bound the nonlinear contribution in terms of the initial distance between the two solutions. Taking $\alpha$ small then preserves the contraction at time $t_c$.
This gives us the existence of a steady state.

\paragraph{Entropy growth.}
We then want to understand the entropy production for this model. We can compute the relative entropy production and find functions for which this is positive when $\alpha>0$ and the stationary state is not a local Maxwellian. We show this is the case when $\tau$ is not a constant.
The derivative of relative entropy is the sum of two nonpositive terms and a term from self collisions that can have either sign. We perturb both the density and the velocity distribution of the stationary state, choosing the perturbation so that the last term is positive and outweighs the other two. The perturbation can be made arbitrarily small, so relative entropy can initially increase even for solutions that start close to the stationary state and later converge back to it.

\paragraph{Plan of the paper.}
Section~\ref{sec:models} introduces the model, and
Section~\ref{sec:boltzmann} develops the linear estimates used in
Section~\ref{sec:nonlinear} to prove global well-posedness and
exponential convergence. Section~\ref{sec:balance-laws} derives balance
laws and Gaussian bounds for the stationary state and shows that it
is not a local Maxwellian when \(\tau\) is nonconstant.
Section~\ref{sec:entropy} proves the results on relative entropy.
Appendix~\ref{app:entropy-growth-proof} justifies the initial growth of relative entropy for mild solutions.
Appendix~\ref{app:numerical-illustrations} describes the figures and the choice of collision rates.

\section{The Thermostatted Boltzmann Model}\label{sec:models}

Equation~\eqref{eq:boltzmann} replaces the BGK relaxation in the model of Evans and Menegaki
\cite{EvansMenegaki2025} by Boltzmann collisions.

\label{sec:coupled-gases}

Both collision operators conserve mass. Self collisions also conserve
the momentum and energy of the gas, while bath collisions can transfer
either to or from the reservoir.

\subsection{The Bath and Collision Law}

For \(u\in\mathbb R^3\) and \(T>0\), let
\[
  M_{u,T}(v)=(2\pi T)^{-3/2}
  \exp\left(-\frac{|v-u|^2}{2T}\right),
  \qquad M_T=M_{0,T}.
\]
Recall the bounds on the reservoir temperature,
\begin{equation}\label{eq:tau-bounds}
  0<\tau_-\le \tau(x)\le \tau_+<\infty
  \qquad (x\in\T).
\end{equation}
We take a representative of \(\tau\) satisfying
\eqref{eq:tau-bounds} at every point.

Both collision operators use the isotropic cutoff kernel
\[
 B_\sigma(v-v_*,\sigma)=\frac{|v-v_*|^\gamma}{4\pi},
 \qquad \sigma\in\mathbb S^2.
\]
Here \(\sigma\) is the outgoing relative direction, and the
post-collision velocities are
\[
 v'=\frac{v+v_*}{2}+\frac{|v-v_*|}{2}\sigma,
 \qquad v_*'=\frac{v+v_*}{2}-\frac{|v-v_*|}{2}\sigma.
\]
The direction \(\sigma\) is uniform on \(\mathbb S^2\).
This gives Maxwell molecules for \(\gamma=0\) and hard spheres for
\(\gamma=1\)
\cite[equations~(1.10)--(1.14)]{BisiCanizoLods2015}.

We also use the
equivalent reflection parametrisation \(n\in\mathbb S^2\):
\[
  v'=v-((v-v_*)\cdot n)n,\qquad
  v_*'=v_*+((v-v_*)\cdot n)n ,
\]
with angular density
\begin{equation}\label{eq:angular-law}
 b_n=\frac{|(v-v_*)\cdot n|}{2\pi|v-v_*|},
 \qquad \int_{\mathbb S^2}b_n\dd n=1,
 \qquad 0\le b_n\le\frac1{2\pi}.
\end{equation}
The outgoing direction and collision normal are related by
\[
  \sigma=\frac{v-v_*-2((v-v_*)\cdot n)n}{|v-v_*|}.
\]

Momentum and energy are conserved,
\(v'+v_*'=v+v_*\) and \(|v'|^2+|v_*'|^2=|v|^2+|v_*|^2\).  For fixed \(x, t\) we write
\(h=h(v)\), \(h_*=h(v_*)\), \(h'=h(v')\), and \(h_*'=h(v_*')\).
For nonnegative \(F,G\) and a test function \(\psi\), the gain operator has
the weak form
\begin{equation}\label{eq:weak-gain-form}
\begin{aligned}
  &\int_{\R}\psi(v)Q^+(F,G)(v)\dd v\\
  &\quad=\int_{\R}\int_{\R}\int_{\mathbb S^2}
  b_n|v-v_*|^\gamma F(v)G(v_*)\psi(v')\dd n\dd v_*\dd v\\
  &\quad=\frac1{4\pi}\int_{\R}\int_{\R}\int_{\mathbb S^2}
  |v-v_*|^\gamma F(v)G(v_*)
  \psi\!\left(\frac{v+v_*}{2}+\frac{|v-v_*|}{2}\sigma\right)
  \dd\sigma\dd v_*\dd v .
\end{aligned}
\end{equation}
We write \(Q(F,G)=Q^+(F,G)-Q^-(F,G)\), where
\[
  Q^-(F,G)(v)=\nu[G](v)F(v),
  \qquad
  \nu[G](v)=\int_{\R}|v-u|^\gamma G(u)\dd u .
\]

We set
\[
 B_n(v-v_*,n)=|v-v_*|^\gamma b_n(v-v_*,n).
\]

Exchanging \(v\) and \(v_*\) in \eqref{eq:weak-gain-form} gives
\(Q^+(F,G)=Q^+(G,F)\). Changing from pre- to post-collision
variables gives
\begin{equation}\label{eq:strong-gain-form}
 Q^+(F,G)(v)=\int_{\R}\int_{\mathbb S^2}
 b_n|v-v_*|^\gamma F(v')G(v_*')\,\dd n\dd v_*.
\end{equation}%
Symmetrising the gain and loss gives the weak forms
\begin{equation}\label{eq:symmetrised-weak-form}
\begin{aligned}
 \int_{\R}\psi\,Q(F,G)\,\dd v
 &=\frac12\int(FG_*-F'G_*')(\psi'-\psi)\,B_n\,\dd n\dd v_*\dd v,\\
 \int_{\R}\psi\,Q(F,F)\,\dd v
 &=\frac14\int(FF_*-F'F_*')(\psi'+\psi_*'-\psi-\psi_*)
   \,B_n\,\dd n\dd v_*\dd v.
\end{aligned}
\end{equation}%
Energy
conservation also gives
\begin{equation}\label{eq:maxwellian-collision-balance}
 M_T(v)M_T(v_*)=M_T(v')M_T(v_*').
\end{equation}

For the kernel estimates in Section~\ref{sec:linear-pointwise}, we use
the Carleman representation \cite[Appendix~C]{GambaPanferovVillani2009}.
For \(T\in[\tau_-,\tau_+]\),
\begin{equation*}
  Q^+(F,M_T)(v)=\int_{\R}k_T(v,u)F(u)\dd u,
\end{equation*}
where, for \(v\ne u\),
\begin{equation}\label{eq:carleman-input}
  k_T(v,u)
  =
  \frac{1}{\pi|v-u|}
  \int_{(v-u)^\perp}
    \bigl(|v-u|^2+|z|^2\bigr)^{(\gamma-1)/2}M_T(v+z)\dd z .
\end{equation}
Here \(u\) is the incoming velocity, \(v\) the outgoing velocity, and
\((v-u)^\perp\) the plane orthogonal to \(v-u\).
We write \(k_x=k_{\tau(x)}\).

\subsection{Types of Solutions}

\begin{definition}[Distributional solution]\label{def:distributional-solution}
A density \(f\) is a distributional solution of
\eqref{eq:boltzmann} on a time interval \(I\) if, for every
\(\varphi\in C_c^1(I^\circ\times\T\times\R)\),
\[
 \int_{I\times\T\times\R}
 f\bigl(\partial_t\varphi+v\cdot\nabla_x\varphi\bigr)\,\dd t\,\dd x\,\dd v
 +\int_{I\times\T\times\R}
 \bigl[(1-\alpha)Q(f,M_{\tau(x)})+\alpha Q(f,f)\bigr]\varphi
 \,\dd t\,\dd x\,\dd v=0.
\]
\end{definition}

The initial condition means \(f(t)\to f_{\rm in}\) in \(L^1\) as
\(t\downarrow0\). Integration along \(x(r)=x-(t-r)v\) gives
\begin{equation}\label{eq:boltzmann-mild}
\begin{aligned}
 f(t,x,v)
 &=\exp\!\left[-\int_{t_0}^t
   \Bigl\{(1-\alpha)\nu[M_{\tau(x-(t-r)v)}](v)
          +\alpha\nu[f](r,x-(t-r)v,v)\Bigr\}\dd r\right]\\
 &\qquad\times
   f(t_0,x-(t-t_0)v,v)\\
 &\quad+\int_{t_0}^t
 \exp\!\left[-\int_s^t
   \Bigl\{(1-\alpha)\nu[M_{\tau(x-(t-r)v)}](v)
          +\alpha\nu[f](r,x-(t-r)v,v)\Bigr\}\dd r\right]\\
 &\qquad\times\Bigl[(1-\alpha)Q^+(f(s),M_{\tau(\cdot)})
   +\alpha Q^+(f(s),f(s))\Bigr](x-(t-s)v,v)\dd s
\end{aligned}
\end{equation}
for almost every \((x,v)\). A mild solution satisfies
\eqref{eq:boltzmann-mild} for every \(t_0\le t\).
Section~\ref{sec:nonlinear} specifies the solution class used here and proves
equivalence with the distributional formulation.

\paragraph{Probabilistic interpretation.}\label{sec:particle-picture}
Between collisions, a particle travels with constant velocity.
At \((x,v)\), its collision rates with bath and gas are
\((1-\alpha)\nu[M_{\tau(x)}](v)\) and \(\alpha\nu[f](t,x,v)\).

The first term in \eqref{eq:boltzmann-mild} represents paths with no collisions after \(t_0\). In the integral, \(s\) is the last collision time. At a bath collision, the new velocity \(v'\) has density
\(k_x(v',v)/\nu[M_{\tau(x)}](v)\) normalised by mass conservation.
At a gas collision, the partner velocity \(v_*\) is drawn from a density proportional to \(|v-v_*|^\gamma f(t,x,v_*)\). We then choose the normal
\(n\) with density \(b_n\) and update the velocities by the collision law.

For \(\alpha=0\), the reservoir fixes collision rates and outgoing velocity distributions, leading to the linear Markov semigroup
\(P_t^\tau\). For \(\alpha>0\), the self-collision rate and the partner velocity distribution depend on the evolving density \(f\), so the evolution is nonlinear.

\subsection{Macroscopic Fields}

We write
\[
  \rho_f=\int_{\R}f\,\dd v,\qquad
  e_f=\int_{\R}|v|^2f\,\dd v.
\]
The momentum density is \(\int_{\R}vf\,\dd v\), and \(e_f/2\) is the
kinetic energy density.  Where \(\rho_f>0\), the temperature is
\[
  T_f=\frac{1}{3\rho_f}
  \left(e_f-\frac1{\rho_f}\left|\int_{\R}vf\,\dd v\right|^2\right).
\]

\subsection{Symmetry Reductions}
\label{sec:symmetry-reductions}

By uniqueness, \(F_\alpha\) inherits symmetries that preserve the
equation and \(\mathcal R\). If \(\tau=\tau(x_1)\), translations in
\(x_2,x_3\) imply that \(F_\alpha\) is independent of these coordinates.
Transport then reduces to \(v_1\partial_{x_1}\), and rotations and
reflections of \((v_2,v_3)\) also preserve the equation and the two
Gaussian norms. Thus \(F_\alpha\) depends only on \(x_1\), \(v_1\), and
\(\sqrt{v_2^2+v_3^2}\).
We use these symmetries in the computations shown in
Figures~\ref{fig:physical-hard-sphere-fields} and~\ref{fig:stationary-velocity}.

\section{Linear Estimates}
\label{sec:boltzmann}
For \(\alpha=0\), the equation is
\[
 \partial_t f+v\cdot\nabla_x f=Q(f,M_{\tau(x)}),
 \qquad f(0)=h.
\]
We write \(P_t^\tau h=f(t)\). We want to control this evolution in
both weighted norms and show that solutions starting from different
probability densities approach one another. The estimates with a source
term in Section~\ref{sec:bath-source} will be used for the nonlinear
equation in Section~\ref{sec:nonlinear}.

\subsection{Construction and Moment Bounds}
\label{sec:linear-construction}

We will test the equation with \(W\) to bound the Gaussian moment
and compare its solution with a multiple of \(W^{-1}\) to obtain
pointwise bounds.

\begin{lemma}[Gaussian moment and pointwise bounds]\label{lem:bath-drift}
For \(W(v)=e^{a|v|^2}\), \(0<a<1/(2\tau_+)\), there are
\(c,C>0\), depending only on \(a,\gamma,\tau_\pm\), such that
\begin{equation}\label{eq:bath-generator-weight}
 \int_{\R}\int_{\mathbb S^2}b_n|v-v_*|^\gamma
 M_{\tau(x)}(v_*)\bigl(W(v')-W(v)\bigr)\dd n\dd v_*
 \le -c\langle v\rangle^\gamma W(v)+C
\end{equation}
and
\begin{equation}\label{eq:bath-gaussian-upper}
 Q(W^{-1},M_{\tau(x)})(v)
 \le\bigl(C-c\langle v\rangle^\gamma\bigr)W^{-1}(v)
\end{equation}
for every \((x,v)\in\T\times\R\).
\end{lemma}

\begin{proof}
In both estimates, a Gaussian factor \(e^{-c(v\cdot n)^2}\)
makes the angular integral concentrate on directions nearly
perpendicular to \(v\). This region narrows as \(|v|\) grows,
giving an extra factor \(1/(1+|v|)\) in the gain. The loss has
size \(\langle v\rangle^\gamma\), so it dominates at large velocities.
For every \(T\in[\tau_-,\tau_+]\),
\begin{equation}\label{eq:hard-frequency-comparison}
 c_{\rm bath}\langle v\rangle^\gamma
 \le\nu[M_T](v)\le C_{\rm bath}\langle v\rangle^\gamma.
\end{equation}
The upper bound follows from the Gaussian moments of \(M_T\).
For the lower bound, restrict to \(|v_*|\le1\) when \(|v|\ge2\),
and to \(3\le|v_*|\le4\) when \(|v|\le2\).
To estimate the gain in \eqref{eq:bath-generator-weight}, use
\[
 \frac{W(v')}{W(v)}
 =e^{-a(v\cdot n)^2+a(v_*\cdot n)^2},
 \qquad
 |v-v_*|^\gamma\le C\langle v\rangle^\gamma\langle v_*\rangle^\gamma.
\]
Since \(a<1/(2\tau_+)\), the integral of
\(\langle v_*\rangle^\gamma M_T(v_*)e^{a|v_*|^2}\) is bounded
uniformly in \(T\). Also,
\[
 \int_{\mathbb S^2}e^{-a(v\cdot n)^2}\dd n
 =2\pi\int_{-1}^1e^{-a|v|^2q^2}\dd q
 \le\frac{C}{1+|v|}.
\]
Together with \(b_n\le1/(2\pi)\), this gives
\begin{equation}\label{eq:bath-weighted-gain-moment}
 \int_{\R}\int_{\mathbb S^2}b_n|v-v_*|^\gamma
 M_{\tau(x)}(v_*)W(v')\dd n\dd v_*
 \le\frac{C\langle v\rangle^\gamma}{1+|v|}W(v).
\end{equation}

For \eqref{eq:bath-gaussian-upper}, the formula
\eqref{eq:strong-gain-form} for \(Q^+\) and energy conservation give
\[
 W(v)W^{-1}(v')M_T(v_*')
 =(2\pi T)^{-3/2}e^{-a|v_*|^2}
   e^{-(1/(2T)-a)|v_*'|^2}.
\]
Here \(|v_*'|^2\ge(v\cdot n)^2\) and
\(1/(2T)-a\ge1/(2\tau_+)-a>0\).
The same angular calculation and integration in \(v_*\) give
\[
 Q^+(W^{-1},M_T)(v)
 \le\frac{C\langle v\rangle^\gamma}{1+|v|}W^{-1}(v).
\]
In each case, subtract the loss in
\eqref{eq:hard-frequency-comparison}.
For large \(|v|\), \(C/(1+|v|)\le c_{\rm bath}/2\), so we may take
\(c=c_{\rm bath}/2\). Increasing \(C\) on the remaining velocity
ball proves both inequalities.
\end{proof}

\begin{proposition}[The linear bath evolution and its moment bounds]
\label{prop:exp-lyapunov}\label{rem:linear-wellposed}
For each \(h\in L^1(W)\), the linear equation has a unique mild
solution in \(C([0,\infty);L^1(W))\) with initial value \(h\).
The solution operators \(P_t^\tau\) form a positive, mass-preserving
semigroup.
With \(c,C\) as in Lemma~\ref{lem:bath-drift}, for every
signed \(h\in L^1(W)\) and \(t\ge0\),
\begin{equation}\label{eq:semigroup-velocity-bound}
  \norm{P_t^\tau h}_{L^1(W)}
  \leq e^{-ct}\norm{h}_{L^1(W)}
  +\frac{C}{c}(1-e^{-ct})\norm{h}_{L^1},
\end{equation}
and
\begin{equation}\label{eq:bath-integrated-moment}
 \|P_t^\tau h\|_{L^1(W)}
 +c\int_0^t\|\langle v\rangle^\gamma P_s^\tau h\|_{L^1(W)}\dd s
 \le \|h\|_{L^1(W)}+Ct\|h\|_{L^1}.
\end{equation}
The semigroup is strongly continuous:
\begin{equation}\label{eq:linear-initial-continuity}
 \lim_{t\downarrow0}\|P_t^\tau h-h\|_{L^1(W)}=0,
 \qquad h\in L^1(W).
\end{equation}
\end{proposition}

\begin{proof}
For \(0\le h\in L^1(W)\), start from \(f^{(0)}=0\) and
iterate \eqref{eq:boltzmann-mild} with \(\alpha=0\), initial datum \(h\)
and \(f^{(n)}\) in the gain. Equivalently, each iterate solves, in the mild sense, \[ (\partial_t+v\cdot\nabla_x)f^{(n+1)} =Q^+(f^{(n)},M_{\tau(x)}) -\nu[M_{\tau(x)}]f^{(n+1)}, \qquad f^{(n+1)}(0)=h. \] Since the gain operator preserves nonnegativity, \(f^{(1)}\ge f^{(0)}=0\). Subtracting two successive mild formulas shows that \(f^{(n+1)}-f^{(n)}\ge0\) whenever \(f^{(n)}-f^{(n-1)}\ge0\). Thus induction gives \(0\le f^{(n)}\le f^{(n+1)}\).

Along \(x=x_0+tv\), write
\(\nu(r)=\nu[M_{\tau(x_0+rv)}](v)\).
Since \(\int_s^L\nu(r)e^{-\int_s^r\nu(q)\dd q}\dd r\le1\),
integrating the mild formula gives, by induction,
\[
 \int_0^L\int_{\T\times\R}\nu[M_\tau]f^{(n+1)}\,\dd v\dd x\dd t
 \le(n+1)\|h\|_{L^1}.
\]
The preceding bound makes both gain and loss integrable, since the velocity integral of the gain equals that of \(\nu[M_\tau]f^{(n)}\). We may therefore integrate against bounded functions of velocity, with the transport term vanishing by periodicity. Taking \(1\) gives \[ \frac{\dd}{\dd t}\int_{\T\times\R}f^{(n+1)} =\int_{\T\times\R}\nu[M_\tau]\bigl(f^{(n)}-f^{(n+1)}\bigr)\le0. \] Hence \(\int f^{(n+1)}(t)\le\|h\|_{L^1}\). To bound the Gaussian moment, we first use \(W_N=\min\{W,N\}\), since the integral with \(W\) is not yet known to be finite. If \(W(v)\ge N\), then \(W_N(v')-W_N(v)\le0\); otherwise, this difference is at most \(W(v')-W(v)\). Using \(f^{(n)}\le f^{(n+1)}\) in the gain and applying \eqref{eq:bath-generator-weight} give
\[
\begin{aligned}
 \frac{\dd}{\dd t}\int_{\T\times\R}W_Nf^{(n+1)}
 &\le\int_{\T\times\R}W_NQ(f^{(n+1)},M_{\tau(x)})\\
 &\le-c\int_{\T\times\R}\langle v\rangle^\gamma
 W{\bf1}_{\{W<N\}}f^{(n+1)}+C\|h\|_{L^1}.
\end{aligned}
\]
Integrate in time and let \(N\to\infty\), then \(n\to\infty\).
Monotone convergence in this inequality and the mild formula gives a
solution \(P_t^\tau h\) satisfying \eqref{eq:bath-integrated-moment}.
By \eqref{eq:bath-weighted-gain-moment} and
\eqref{eq:hard-frequency-comparison}, the gain and loss belong to
\(L^1((0,L);L^1(W))\). Testing the limit equation with \(1\) proves
mass conservation. Testing with \(W\) is now justified and gives
\[
 \frac{\dd}{\dd t}\int_{\T\times\R}Wf
 +c\int_{\T\times\R}\langle v\rangle^\gamma Wf
 \le C\int_{\T\times\R}f.
\]
Since \(\langle v\rangle^\gamma\ge1\) and \(\int f(t)=\int h\),
Gr\"onwall's inequality gives \eqref{eq:semigroup-velocity-bound}.

Writing \((T_th)(x,v)=h(x-tv,v)\), the integral equation also gives
\[
 \|P_t^\tau h-T_th\|_{L^1(W)}
 \le C\int_0^t
 \|\langle v\rangle^\gamma P_s^\tau h\|_{L^1(W)}\dd s
 \longrightarrow0\qquad(t\downarrow0).
\]
Continuity of translations proves
\eqref{eq:linear-initial-continuity} for nonnegative data.

The equation also gives \(L^1\) contraction.
For the difference \(w\) of two solutions bounded in \(L^1(W)\),
taking absolute values and integrating gives
\[
 \|w(t)\|_{L^1}
 \le\|w(0)\|_{L^1}
 +\int_0^t\int_{\T\times\R}Q(|w(s)|,M_{\tau(x)})\,\dd v\dd x\dd s
 =\|w(0)\|_{L^1}.
\]
The bound \(\nu[M_\tau]\le CW\) and the \(L^1(W)\) bound on \(w\)
give \(Q^\pm(|w|,M_{\tau(\cdot)})\in L^1((0,t)\times\T\times\R)\)
for every \(t\). By linearity, the difference of two solutions with the same initial datum solves the equation with zero initial datum. The estimate above forces this difference to vanish.
Extend the construction linearly to signed data. Positivity gives
\(|P_t^\tau h|\le P_t^\tau|h|\), hence both moment bounds.
Continuity follows by treating the positive and negative parts of
\(h\) separately.

\end{proof}

\begin{remark}
Positivity gives \(|P_t^\tau h|\le P_t^\tau|h|\), so mass conservation implies
\[
 \|P_t^\tau h\|_{L^1}\le\int_{\T\times\R}P_t^\tau|h|\,\dd v\dd x
 =\|h\|_{L^1},\qquad h\in L^1(W).
\]
By linearity, the same estimate applies to differences of initial data.
Since \(L^1(W)\) is dense in \(L^1\), we can therefore define
\(P_t^\tau h\) for \(h\in L^1\) by approximation in \(L^1\).
\end{remark}

\subsection{Weighted \texorpdfstring{\(L^1\)}{L1} Contraction by Harris's Theorem}
\label{sec:linear-harris}

Together with \eqref{eq:semigroup-velocity-bound}, a lower bound on
transition probabilities from bounded initial velocities will allow us to
apply Harris's theorem
\cite[Theorem~1.3 and Remark~1.4]{HairerMattingly2011}.
The same iteration works for a Dirac initial measure, so we define
\(p_t((x,v),\cdot)=P_t^\tau\delta_{(x,v)}\).

Starting from \((x_0,v_0)\), the first collision produces a positive
density on a ball of velocities, while the position is still fixed by the
starting point and collision time. During the following flight, a large
enough ball of velocities reaches the whole torus. Another collision then
gives a velocity density at each position. We follow the construction in
\cite[Lemmas~4.2--4.3]{CanizoCaoEvansYoldas2020}.

\begin{lemma}\label{lem:local-gain-lower}
For every \(R>0\), there is \(c_R>0\) such that
\begin{equation}\label{eq:bath-kernel-lower}
 k_{\tau(x)}(v,u)
 \ge c_R\langle v\rangle^{\gamma-2}M_{\tau(x)}(v),
 \qquad |u|\le R,\quad v\ne u,
\end{equation}
for every \(x\). Consequently, for every \(r>0\), there is
\(c>0\) such that every finite nonnegative measure \(\mu\) on velocity
space satisfies
\[
  Q^+(\mu,M_{\tau(x)})(\dd v)
  \ge c\,\mu(\{|u|\le R\})\one_{\{|v|\le r\}}\dd v
\]
for every \(x\).
\end{lemma}

\begin{proof}
We use a lower-bound argument based on the Carleman representation
\cite[Lemma~4.2]{CanizoCaoEvansYoldas2020}. In
\eqref{eq:carleman-input}, restrict the plane integral to
\(1\le |z|\le2\).
For \(|u|\le R\), put \(d=|v-u|>0\). On this annulus,
\(v\cdot z=u\cdot z\), so
\(M_{\tau(x)}(v+z)\ge c_RM_{\tau(x)}(v)\). Moreover,
\[
 d\le C_R\langle v\rangle,\qquad
 C_R^{-1}\langle v\rangle^2
 \le d^2+|z|^2
 \le C_R\langle v\rangle^2.
\]
Equation~\eqref{eq:carleman-input} now gives
\eqref{eq:bath-kernel-lower}. Its right side is bounded below
by a positive constant on \(|v|\le r\). Integrate in
\(\mu(\dd u)\) over \(|u|\le R\) to obtain the claim. 
\end{proof}

On measures, \(T_t\) acts by
\((x,v)\mapsto(x+tv,v)\).
\begin{lemma}\label{lem:transport-cover}
Let \(t_0>0\) and \(r\geq\sqrt3/(2t_0)\).
For every \(t\geq t_0\) and \(x_0\in\T\), the spatial marginal of
\[
  T_t\bigl(\delta_{x_0}(\dd x)\one_{\{|w|\le r\}}\dd w\bigr)
\]
dominates \(t^{-3}\dd x\).
\end{lemma}

\begin{proof}
Periodising \(x=x_0+tw\) gives the torus covering calculation from
\cite[Lemma~3.2]{CanizoCaoEvansYoldas2020}:
\[
  T_t\bigl(\delta_{x_0}(\dd x)\one_{\{|w|\le r\}}\dd w\bigr)
  =
  t^{-3}\sum_{k\in\Z}
  \one_{\{|w_k(x)|\le r\}}\dd x\,\delta_{w_k(x)}(\dd w),
  \qquad w_k(x)=\frac{x-x_0-k}{t}.
\]
Every \(x\) has a representative with
\(|x-x_0-k|\leq\sqrt3/2\), so \(|w_k(x)|\leq\sqrt3/(2t)\leq r\).
The sum therefore contains at least one term at each \(x\).
\end{proof}

\begin{proposition}[Uniform lower bound on transition probabilities]\label{prop:linear-minorisation}
Fix \(t_\ast>0\) and \(R_0<\infty\).  There are
\(\varepsilon>0\) and \(r>0\)
such that
\[
  p_{t_\ast}((x_0,v_0),\dd x\dd v)
  \geq\varepsilon\,\one_{\{|v|\leq r\}}\dd x\dd v
\]
for all \(x_0\in\T\) and \(|v_0|\leq R_0\).  This gives a
uniform lower bound on every set on which \(W\) is bounded.
\end{proposition}

\begin{proof}
Take \(\mu=\delta_{(x_0,v_0)}\), with \(|v_0|\leq R_0\).
Put \(t_0=t_\ast/3\), take \(r\) as in
Lemma~\ref{lem:transport-cover}, and choose
\(R\geq\max\{R_0,r\}\). Set
\[
  \mathcal Gg
  =\one_{\{|v|\leq r\}}
  Q^+\bigl(\one_{\{|u|\leq R\}}g,M_{\tau(x)}\bigr).
\]
All three free flights have velocities in \(\{|v|\le R\}\),
where \eqref{eq:hard-frequency-comparison} bounds the loss rate.
Keeping two gains, both restricted to \(\mathcal G\), and bounding
the total loss factor from below gives
\[
  P_{t_\ast}^\tau \mu
  \geq e^{-Ct_\ast}
  \int_0^{t_\ast}\int_0^s
  T_{t_\ast-s}\mathcal G T_{s-\ell}\mathcal G T_\ell \mu
  \,\dd \ell\dd s
\]
in the sense of measures.

The first application of \(\mathcal G\) gives, for some \(c>0\),
\[
  \mathcal G T_\ell \mu
  \geq c\,\delta_{x_0+\ell v_0}(\dd x)
  \one_{\{|v|\leq r\}}\dd v .
\]
For \(s-\ell\geq t_0\), Lemma~\ref{lem:transport-cover} and a second
application of \(\mathcal G\) give
\[
  \mathcal G T_{s-\ell}\mathcal G T_\ell \mu
  \geq c^2(s-\ell)^{-3}\dd x\,
  \one_{\{|v|\leq r\}}\dd v .
\]
Restrict the time integrals to \(0<\ell<t_0\) and
\(2t_0<s<t_\ast\), so \(t_0\le s-\ell\le t_\ast\).
The final transport preserves
\(\dd x\,\one_{\{|v|\le r\}}\dd v\).
Integration over this rectangle proves the lower bound.
\end{proof}

\begin{proposition}[Contraction for the linear evolution]
\label{prop:bath-harris}
There are \(C,\lambda>0\), depending only on \(a\), \(\gamma\),
\(\tau_\pm\) and the torus, such that, for every \(t\ge0\) and
signed \(h\in L^1(W)\) with \(\int_{\T\times\R}h=0\),
\begin{equation}\label{eq:hard-harris-contraction}
  \norm{P_t^\tau h}_{L^1(W)}
  \leq Ce^{-\lambda t}\norm{h}_{L^1(W)}.
\end{equation}
\end{proposition}

\begin{proof}
Fix \(t_H>0\). Applying \eqref{eq:semigroup-velocity-bound}
to \(\delta_{(x,v)}\) gives the drift condition
\[
 \int_{\T\times\R}W(u)\,p_{t_H}((x,v),\dd y\dd u)
 \le e^{-ct_H}W(v)+\frac Cc(1-e^{-ct_H}).
\]
Proposition~\ref{prop:linear-minorisation} gives minorisation
on every sufficiently large set \(\{W\le R\}\).
Harris's theorem \cite{Harris1956}, in the form
\cite[Theorem~1.3]{HairerMattingly2011}, therefore gives
\(\beta>0\) and \(q\in(0,1)\) such that
\[
 \|P_{t_H}^\tau h\|_{L^1(1+\beta W)}
 \le q\|h\|_{L^1(1+\beta W)},
 \qquad \int_{\T\times\R}h=0.
\]
Since \(W\ge1\), the weights \(W\) and \(1+\beta W\) give
equivalent norms. Iterating this inequality gives exponential decay
at multiples of \(t_H\); \eqref{eq:semigroup-velocity-bound} controls
the remaining time interval. This proves
\eqref{eq:hard-harris-contraction}.
\end{proof}

\subsection{Weighted \texorpdfstring{\(L^\infty\)}{L-infinity} Bounds and Contraction}
\label{sec:linear-pointwise}

By \eqref{eq:hard-frequency-comparison}, the exponential in
Duhamel's formula is at most
\(e^{-c\langle v\rangle^\gamma(t-s)}\), and
\[
 \int_0^t e^{-c\langle v\rangle^\gamma(t-s)}
             \langle v\rangle^\gamma\dd s\le\frac1c.
\]
Thus time integration compensates for the factor
\(\langle v\rangle^\gamma\) in the gain estimate. To control the
velocity integral, we use the following bounds from the Carleman
representation \eqref{eq:carleman-input}.

\begin{lemma}[Bounds for the bath kernel with a Gaussian weight]\label{lem:weighted-carleman}
Assume \(0<a<1/(2\tau_+)\).  There are constants \(C,c>0\), depending only on
\(a\), \(\tau_\pm\), and \(\gamma\), such that
\[
  \frac{e^{a(|v|^2-|u|^2)}}{\langle v\rangle^\gamma}
  k_T(v,u)\le\widetilde k(v,u),
  \qquad T\in[\tau_-,\tau_+],
\]
where, for \(v\ne u\),
\begin{equation}\label{eq:weighted-kernel-definition}
  \widetilde k(v,u)
  =
  \frac{C}{|v-u|^2}
  \exp\!\left[-c\left(|v-u|^2+
  \frac{(|v|^2-|u|^2)^2}{|v-u|^2}\right)\right],
\end{equation}
and \(\widetilde k(v,v)=0\).  Moreover,
\begin{equation}\label{eq:K-bounds}
  \sup_v(1+|v|)\int_{\R}\widetilde k(v,u)\dd u\le C,
  \qquad
  \sup_u(1+|u|)\int_{\R}\widetilde k(v,u)\dd v\le C,
\end{equation}
and, for every \(\varepsilon>0\),
\begin{equation}\label{eq:K-local}
  \one_{\{|v-u|\ge\varepsilon\}}\widetilde k(v,u)
  \le C_\varepsilon e^{-c|v-u|^2},
  \qquad
  \sup_v\int_{\{|u-v|<\varepsilon\}}\widetilde k(v,u)\dd u
  \le C\varepsilon .
\end{equation}
Consequently, for every nonnegative velocity function \(h\), uniformly in
\(x\) and for almost every \(v\),
\begin{equation}\label{eq:weighted-bath-gain}
  W(v)\langle v\rangle^{-\gamma}
  Q^+(W^{-1}h,M_{\tau(x)})(v)
  \le \int_{\R}\widetilde k(v,u)h(u)\dd u .
\end{equation}
\end{lemma}

\begin{proof}
Put \(r=|v-u|\), \(\Delta=|v|^2-|u|^2\), and \(\beta=1/(8T)\).
If \(v_\perp\) is the component of \(v\) orthogonal to \(v-u\), then
\[
 |v+z|^2=\frac14\left(r+\frac{\Delta}{r}\right)^2
          +|v_\perp+z|^2,
 \qquad z\perp(v-u).
\]
For every \(\gamma\ge0\),
\(r^{-1}(r^2+|z|^2)^{(\gamma-1)/2}
 \le r^{-2}(r^2+|z|^2)^{\gamma/2}\).
After the translation \(y=v_\perp+z\), the plane integral satisfies
\[
 \int_{(v-u)^\perp}
 \frac{(r^2+|y-v_\perp|^2)^{\gamma/2}}{\langle v\rangle^\gamma}
 e^{-|y|^2/(2T)}\dd y
 \le C\langle r\rangle^\gamma
       \int_{(v-u)^\perp}\langle y\rangle^\gamma e^{-|y|^2/(2T)}\dd y
 \le C\langle r\rangle^\gamma.
\]
Consequently, \eqref{eq:carleman-input} gives
\[
 \frac{e^{a\Delta}}{\langle v\rangle^\gamma}k_T(v,u)
 \le\frac{C\langle r\rangle^\gamma}{r^2}
       e^{a\Delta-\beta(r+\Delta/r)^2}.
\]
The identity
\[
 \beta\left(r+\frac{\Delta}{r}\right)^2-a\Delta
 =\left(\beta-\frac a4\right)\left(r+\frac{\Delta}{r}\right)^2
   +\frac a4\left(r-\frac{\Delta}{r}\right)^2
 \ge c\left(r^2+\frac{\Delta^2}{r^2}\right)
\]
holds uniformly in \(T\in[\tau_-,\tau_+]\), because
\(0<a<1/(2\tau_+)\). Absorb \(\langle r\rangle^\gamma\) into the
Gaussian to obtain the kernel bound.

For \eqref{eq:K-bounds}, fix \(u\) and write \(v=u+r\omega\).
Then \(\Delta/r=2u\cdot\omega+r\), and
\[
 \int_{\mathbb S^2}e^{-c(2u\cdot\omega+r)^2}\dd\omega
 =2\pi\int_{-1}^1e^{-c(2|u|q+r)^2}\dd q
 \le\frac{C}{1+|u|}.
\]
This integral is bounded both by \(4\pi\) and, for \(u\ne0\),
by \(C/|u|\). Polar coordinates therefore give
\[
 \int_{\R}\widetilde k(v,u)\dd v
 \le\frac{C}{1+|u|}\int_0^\infty e^{-cr^2}\dd r
 \le\frac{C}{1+|u|}.
\]
By symmetry, the same bound holds with \(u\) and \(v\) exchanged.
For \eqref{eq:K-local}, bound \(r^{-2}\) by \(\varepsilon^{-2}\)
outside the ball and integrate \(r^{-2}r^2\) from zero to
\(\varepsilon\) inside it. Finally, the gain representation with kernel
\eqref{eq:carleman-input}, applied to \(W^{-1}h\) for \(h\ge0\), gives
\eqref{eq:weighted-bath-gain}.
\end{proof}

\begin{lemma}\label{lem:bath-Linf-boundedness}
For every \(L>0\), there is \(C_L<\infty\) such that
\[
 \sup_{0\le t\le L}\|P_t^\tau h\|_{L^\infty(W)}
 \le C_L\|h\|_{L^\infty(W)}
\]
for every signed \(h\in L^\infty(W)\).
\end{lemma}
\begin{proof}
By \eqref{eq:bath-gaussian-upper},
\(e^{Ct}\|h\|_{L^\infty(W)}W^{-1}\) is a supersolution and bounds
\(|h|\) initially. Integration along characteristics and positivity of the gain
bound every Duhamel iterate by this Gaussian.
Apply the comparison first to
\(|h|\one_{\{|v|\le N\}}\in L^1(W)\). The Gaussian bound on \(|h|\)
makes these truncations converge in \(L^1\). By \(L^1\) contraction,
their evolutions converge to \(P_t^\tau|h|\) in \(L^1\). A subsequence
converges almost everywhere, so the Gaussian bound passes to the limit. Since \(|P_t^\tau h|\le P_t^\tau|h|\),
we obtain the result with \(C_L=e^{CL}\).
\end{proof}

Applying Duhamel's formula once gives a velocity integral at a position independent of the integration velocity. Substituting the formula again makes the earlier position depend on the intermediate velocity, allowing us to change variables and use the
\(L^1(W)\) bound. This argument can be found, for example, in Guo \cite{Guo2010}
and Vidav \cite{Vidav1970}. For differences of probability densities, the
\(L^1(W)\) decay from Harris's theorem then gives contraction in the sum
of the two norms.

\begin{proposition}[Estimate in \(L^\infty(W)\)]
\label{prop:pointwise-recovery}
\label{cor:linear-pointwise-return}
\label{cor:linear-intersection-contraction}
There is \(C<\infty\) such that
\begin{equation}\label{eq:linear-Linf-uniform}
 \sup_{t\ge0}\norm{P_t^\tau h}_{L^\infty(W)}
 \le C\bigl(\norm{h}_{L^1(W)}+\norm{h}_{L^\infty(W)}\bigr).
\end{equation}
Moreover, for every \(\varepsilon>0\), there are
\(L_0,C<\infty\), depending on \(\varepsilon\), such that,
for \(L\ge L_0\),
\begin{equation}\label{eq:linear-Linf-return}
 \norm{P_L^\tau h}_{L^\infty(W)}
 \le \varepsilon\norm{h}_{L^\infty(W)}
 +C\norm{h}_{L^1(W)}.
\end{equation}
Both estimates hold for signed
\(h\in L^1(W)\cap L^\infty(W)\).

For every \(\theta\in(0,1)\), there is
\(L_{\theta}>0\) such that, for every
\(L\ge L_{\theta}\) and every
\(h\in L^1(W)\cap L^\infty(W)\) with \(\int_{\T\times\R}h=0\),
\begin{equation}\label{eq:linear-intersection-contraction}
\begin{aligned}
 \|P_L^\tau h\|_{L^1(W)}+\|P_L^\tau h\|_{L^\infty(W)}
 \le\theta\bigl(\|h\|_{L^1(W)}+\|h\|_{L^\infty(W)}\bigr).
\end{aligned}
\end{equation}
\end{proposition}
\begin{proof}
Let \(h_0\ge0\) belong to \(L^1(W)\cap L^\infty(W)\), and set
\(h(t)=P_{t-t_0}^\tau h_0\). For every \(\varepsilon>0\), we first
prove the following estimate for all \(t\ge t_0\):
\begin{equation}\label{eq:pointwise-recovery}
\begin{aligned}
 \norm{h(t)}_{L^\infty(W)}
 &\le Ce^{-c(t-t_0)}\norm{h_0}_{L^\infty(W)}\\
 &\quad+C\sup_{t_0\le q\le t}\norm{h(q)}_{L^1(W)}
 +\varepsilon\sup_{t_0\le q\le t}\norm{h(q)}_{L^\infty(W)}.
\end{aligned}
\end{equation}
Here \(C\) may depend on \(\varepsilon\), but the constants
are independent of \(t_0,t\), and \(h_0\).
Fix \(t\ge t_0\) and write \(J=[t_0,t]\).

For \(t\ge s\), set
\[
 E(t,s;x,v)
 =\exp\!\left[-\int_s^t
 \nu[M_{\tau(x-(t-r)v)}](v)\dd r\right].
\]
By \eqref{eq:hard-frequency-comparison},
\begin{equation}\label{eq:loss-factor-bound}
 0\le E(t,s;x,v)
 \le e^{-c\langle v\rangle^\gamma(t-s)}.
\end{equation}
Applying Duhamel's formula twice gives
\[
 h(t,x,v)=I_0(t,x,v)+I_1(t,x,v)+I_2(t,x,v),
\]
where
\[
 I_0(t,x,v)
 =E(t,t_0;x,v)h_0(x-(t-t_0)v,v),
\]
\begin{align*}
 I_1(t,x,v)
 &=\int_{t_0}^t\int_{\R}
 E(t,s;x,v)k_{x-(t-s)v}(v,u)\,
 E(s,t_0;x-(t-s)v,u)\\
 &\qquad\times
 h_0(x-(t-s)v-(s-t_0)u,u)\dd u\dd s.
\end{align*}
The remaining integral is \(I_2\ge0\); it contains
\(h(r)\) at an earlier time rather than \(h_0\).
As in the probabilistic interpretation in Section~\ref{sec:particle-picture}, the three terms describe trajectories with zero, one, and at least
two collisions between \(t_0\) and \(t\).

Set \(H=\norm{h_0}_{L^\infty(W)}\).  From
\eqref{eq:loss-factor-bound},
\[
 W(v)I_0(t,x,v)\le e^{-c(t-t_0)}H.
\]
For \(I_1\), \eqref{eq:loss-factor-bound},
\eqref{eq:weighted-bath-gain}, and the first bound in \eqref{eq:K-bounds} give
\[
\begin{aligned}
 W(v)I_1(t,x,v)
 &\le CH\int_{t_0}^t
 e^{-c\langle v\rangle^\gamma(t-s)}\langle v\rangle^\gamma
 \int_{\R}\widetilde k(v,u)
 e^{-c\langle u\rangle^\gamma(s-t_0)}\dd u\dd s\\
 &\le CH e^{-c(t-t_0)/2}\int_{t_0}^t
 e^{-c\langle v\rangle^\gamma(t-s)/2}
 \langle v\rangle^\gamma\dd s
 \le CH e^{-c(t-t_0)/2}.
\end{aligned}
\]
Here \(\langle v\rangle^\gamma,\langle u\rangle^\gamma\ge1\)
give the factor \(e^{-c(t-t_0)/2}\), while the remaining time
integral is at most \(2/c\).

Applying \eqref{eq:weighted-bath-gain} twice, and using
\eqref{eq:loss-factor-bound}, gives
\[
\begin{aligned}
 W(v)I_2(t,x,v)
 &\le C\int_{t_0}^t\int_{t_0}^s\int_{\R\times\R}
 e^{-c\{\langle v\rangle^\gamma(t-s)
         +\langle u\rangle^\gamma(s-r)\}}
 \langle v\rangle^\gamma\langle u\rangle^\gamma\\
 &\quad\times\widetilde k(v,u)\widetilde k(u,w)
 W(w)h\bigl(r,x-(t-s)v-(s-r)u,w\bigr)
 \,\dd w\dd u\dd r\dd s .
\end{aligned}
\]

To use the \(L^1(W)\) bound in \(I_2\), we need to integrate \(h\)
in position as well as velocity. Changing \(u\) to the position at time
\(r\) does this, but introduces a factor \((s-r)^{-3}\).
We first use the supremum norm to estimate short time gaps, large
velocities and small differences between successive velocities.

Fix \(R>1\) and \(0<\delta<1\). First consider the region
\[
 |v|>R,\qquad |u|>2R,\qquad\text{or}\qquad |w|>3R.
\]
Bound \(W(w)h(r,\ldots,w)\) by
\(\sup_{q\in J}\|h(q)\|_{L^\infty(W)}\).
The two time integrals satisfy
\[
 \int_{t_0}^t e^{-c\langle v\rangle^\gamma(t-s)}
                  \langle v\rangle^\gamma
 \int_{t_0}^s e^{-c\langle u\rangle^\gamma(s-r)}
                  \langle u\rangle^\gamma\dd r\dd s
 \le\frac1{c^2}.
\]
If \(|v|>R\), the first
estimate in \eqref{eq:K-bounds} gives a factor \(C/R\).  If
\(|v|\le R\) and \(|u|>2R\), then \(|v-u|>R\); if \(|u|\le2R\) and
\(|w|>3R\), then \(|u-w|>R\).  In the latter two cases the Gaussian bound in
\eqref{eq:K-local} gives \(Ce^{-cR^2}\).  Hence this region contributes at
most
\[
 \left(\frac C{R}+Ce^{-cR^2}\right)
 \sup_{q\in J}\norm{h(q)}_{L^\infty(W)}.
\]

We now restrict to \(|v|\le R\), \(|u|\le2R\) and
\(|w|\le3R\). For the part with \(s-r\le\delta\),
\[
 \int_{t_0}^s
 e^{-c\langle u\rangle^\gamma(s-r)}\langle u\rangle^\gamma
 \one_{\{s-r\le\delta\}}\dd r
 \le \langle2R\rangle^\gamma\delta.
\]
For the parts with \(|v-u|<\delta\) or
\(|u-w|<\delta\), \eqref{eq:K-local} gives
\[
 \sup_v\int_{|v-u|<\delta}\widetilde k(v,u)\dd u
 +\sup_u\int_{|u-w|<\delta}\widetilde k(u,w)\dd w
 \le C\delta.
\]
Together these regions contribute at most
\[
 C_R\delta\sup_{q\in J}\norm{h(q)}_{L^\infty(W)}.
\]

On the remaining region,
\[
\begin{gathered}
 |v|\le R,\qquad |u|\le2R,\qquad |w|\le3R,\\
 |v-u|,|u-w|\ge\delta,\qquad s-r\ge\delta.
\end{gathered}
\]
Put \(\ell=s-r\) and set
\(y=x-(t-s)v-\ell u\), so that \(\dd u=\ell^{-3}\dd y\).
Since \(y\) is a torus variable, the velocities giving the same
position are
\[
 u_k=\frac{x-(t-s)v-y-k}{\ell},\qquad k\in\mathbb Z^3.
\]
In the sum below retain only the indices satisfying
\(|u_k|\le2R\) and \(|v-u_k|,|u_k-w|\ge\delta\).
There are at most \((4R\ell+2)^3\) such indices, and
\eqref{eq:K-local} bounds the product of kernels by \(C_\delta\).
Hence
\[
\begin{aligned}
 \sum_k\ell^{-3}\widetilde k(v,u_k)\widetilde k(u_k,w)
 &\le C_\delta\ell^{-3}(4R\ell+2)^3\\
 &\le C_\delta(R^3+\ell^{-3})\le C_{\delta,R}.
\end{aligned}
\]
This bound is uniform in the time gap \(\ell\ge\delta\).

On this region, \(\langle u\rangle^\gamma\le\langle2R\rangle^\gamma\) and
\(e^{-c\langle u\rangle^\gamma(s-r)}\le e^{-c(s-r)}\).
The integrals in \(y\) and \(w\) are now controlled by
\(\|h(r)\|_{L^1(W)}\), so this part of \(I_2\) is bounded by
\begin{align*}
 &C_{\delta,R}\sup_{q\in J}\norm{h(q)}_{L^1(W)}
 \int_{t_0}^t e^{-c\langle v\rangle^\gamma(t-s)}
 \langle v\rangle^\gamma
 \int_{t_0}^s e^{-c(s-r)}\dd r\dd s\\
 &\qquad\le
 C_{\delta,R}\sup_{q\in J}\norm{h(q)}_{L^1(W)}.
\end{align*}
Combining the three estimates gives
\[
 \norm{I_2(t)}_{L^\infty(W)}
 \le\left(\frac C{R}+Ce^{-cR^2}+C_R\delta\right)
 \sup_{q\in J}\norm{h(q)}_{L^\infty(W)}
 +C_{\delta,R}\sup_{q\in J}\norm{h(q)}_{L^1(W)}.
\]
Choose \(R\) so that
\(C/R+Ce^{-cR^2}\le\varepsilon/2\), then choose \(\delta\) so that
\(C_R\delta\le\varepsilon/2\).
Together with the bounds for \(I_0\) and \(I_1\), this proves
\eqref{eq:pointwise-recovery}.

For nonnegative initial data \(h\), the moment bound
\eqref{eq:semigroup-velocity-bound} and \(W\ge1\) give
\begin{equation}\label{eq:linear-L1-uniform}
 \sup_{t\ge0}\|P_t^\tau h\|_{L^1(W)}
 \le C_1\|h\|_{L^1(W)}.
\end{equation}
On \([0,L]\), Lemma~\ref{lem:bath-Linf-boundedness} makes the
\(L^\infty(W)\) supremum finite.
Take the supremum over \(0\le t\le L\) in
\eqref{eq:pointwise-recovery}, with \(t_0=0\) and its final
coefficient equal to \(1/2\). Using \eqref{eq:linear-L1-uniform}
and moving the last term to the left gives, uniformly in \(L\),
\[
 \sup_{0\le t\le L}\|P_t^\tau h\|_{L^\infty(W)}
 \le C_U\bigl(\|h\|_{L^1(W)}+\|h\|_{L^\infty(W)}\bigr).
\]
Letting \(L\to\infty\) proves \eqref{eq:linear-Linf-uniform}.

Now fix \(\varepsilon>0\) and apply
\eqref{eq:pointwise-recovery} at \(t_0=0\), \(t=L\), with its
final coefficient chosen as \(\varepsilon/(2C_U)\).
Choose \(L_0\) so large that \(Ce^{-cL}\le\varepsilon/2\)
whenever \(L\ge L_0\). Equations \eqref{eq:linear-L1-uniform}
and \eqref{eq:linear-Linf-uniform} then give
\[
 \|P_L^\tau h\|_{L^\infty(W)}
 \le\varepsilon\|h\|_{L^\infty(W)}+C\|h\|_{L^1(W)},
\]
which proves \eqref{eq:linear-Linf-return}.  For signed \(h\), use
\(|P_t^\tau h|\le P_t^\tau|h|\).

For \(\int h=0\), choose \(\ell>0\) so that
\eqref{eq:linear-Linf-return} holds with
\(\varepsilon=\theta/(2C_U)\).
For \(L\ge\ell\), apply it on the final interval of length \(\ell\):
\[
\begin{aligned}
 \norm{P_L^\tau h}_{L^\infty(W)}
 &\le\frac{\theta}{2C_U}
       \norm{P_{L-\ell}^\tau h}_{L^\infty(W)}
       +C\norm{P_{L-\ell}^\tau h}_{L^1(W)}\\
 &\le\frac{\theta}{2}
       \bigl(\norm{h}_{L^1(W)}+\norm{h}_{L^\infty(W)}\bigr)
       +Ce^{-\lambda(L-\ell)}\norm{h}_{L^1(W)}.
\end{aligned}
\]
Here we used \eqref{eq:linear-Linf-uniform} for the first
term and Harris's theorem for the second. The latter also gives
\[
 \norm{P_L^\tau h}_{L^1(W)}
 \le Ce^{-\lambda L}\norm{h}_{L^1(W)}.
\]
Choose \(L_\theta\ge\ell\) so that the two exponentially
decaying coefficients sum to at most \(\theta/2\) for
\(L\ge L_\theta\). Adding the bounds proves
\eqref{eq:linear-intersection-contraction}.
\end{proof}
\subsection{Bounds for the Linear Equation with a Source}
\label{sec:bath-source}

\begin{lemma}[Bounds for the equation with a source]
\label{lem:bath-source-duhamel}
Let \(L>0\) and let \(\Phi\) be a measurable signed source. Set
\[
 u(t)=\int_0^tP_{t-s}^\tau\Phi(s)\dd s,
 \qquad 0\le t\le L.
\]
Then
\begin{equation}\label{eq:bath-source-duhamel}
 \sup_{0\le t\le L}\|u(t)\|_{L^\infty(W)}
 \le C\operatorname*{ess\,sup}_{0\le s\le L}
 \|\langle v\rangle^{-\gamma}\Phi(s)\|_{L^\infty(W)}.
\end{equation}
Also,
\begin{align}
 &\sup_{0\le t\le L}\|u(t)\|_{L^1(W)}
 +\int_0^L
   \|\langle v\rangle^\gamma u(s)\|_{L^1(W)}\dd s
 \notag\\
 &\qquad\le C\int_0^L\|\Phi(s)\|_{L^1(W)}\dd s.
 \label{eq:bath-source-one-integrated}
\end{align}
Each estimate holds whenever its right-hand side is finite.
The constant \(C\) may depend on \(L\), but not on \(\Phi\).
Under either assumption, \(u\) is a mild solution of
\[
 \partial_tu+v\cdot\nabla_xu=Q(u,M_{\tau(x)})+\Phi,
 \qquad u(0)=0.
\]
If \(\Phi\in L^1((0,L);L^1(W))\), then
\(u\in C([0,L];L^1(W))\). In that case, it is the unique distributional
solution bounded in \(L^1(W)\) whose initial value is zero in \(L^1\).
\end{lemma}

\begin{proof}
Either assumption gives
\(\Phi\in L^1((0,L)\times\T\times\R)\).
For \(\Phi\ge0\), substituting \eqref{eq:boltzmann-mild} with
\(\alpha=0\) into the definition of \(u\) and applying Tonelli's
theorem gives the mild formula. Apply this to the positive and
negative parts for signed \(\Phi\).
Positivity of \(P_t^\tau\) reduces both estimates to \(\Phi\ge0\).

For \eqref{eq:bath-source-duhamel}, linearity lets us assume
without loss of generality that
\(0\le\Phi\le\langle v\rangle^\gamma W^{-1}\).
The pointwise Gaussian inequality \eqref{eq:bath-gaussian-upper}
gives
\[
 (\partial_t+v\cdot\nabla_x)(c^{-1}e^{Ct}W^{-1})
 -Q(c^{-1}e^{Ct}W^{-1},M_{\tau(x)})
 \ge e^{Ct}\langle v\rangle^\gamma W^{-1}\ge\Phi.
\]
Comparison with this Gaussian, as in
Lemma~\ref{lem:bath-Linf-boundedness}, proves
\eqref{eq:bath-source-duhamel}.

For \eqref{eq:bath-source-one-integrated}, assume
\(\int_0^L\|\Phi(r)\|_{L^1(W)}\dd r<\infty\).
The bound \eqref{eq:semigroup-velocity-bound} gives
\[
 \sup_{0\le t\le L}\|u(t)\|_{L^1(W)}
 \le C\int_0^L\|\Phi(r)\|_{L^1(W)}\dd r.
\]
For the time integral, Tonelli's theorem and
\eqref{eq:bath-integrated-moment} give
\[
\begin{aligned}
 \int_0^L\|\langle v\rangle^\gamma u(s)\|_{L^1(W)}\dd s
 &\le\int_0^L\int_r^L
 \|\langle v\rangle^\gamma P_{s-r}^\tau\Phi(r)\|_{L^1(W)}\dd s\dd r\\
 &\le C\int_0^L\|\Phi(r)\|_{L^1(W)}\dd r.
\end{aligned}
\]
Together these prove \eqref{eq:bath-source-one-integrated}.

Strong continuity of \(P_t^\tau\) on \(L^1(W)\) gives
\(u\in C([0,L];L^1(W))\).
Uniqueness follows by applying Proposition~\ref{prop:exp-lyapunov}
to the difference of two solutions.
\end{proof}

\section{Nonlinear Well-Posedness and Relaxation to the Stationary State}
\label{sec:nonlinear}
\label{sec:nonlinear-closure}

With the collision partner's distribution fixed, we solve the linear equation by iterating Duhamel's formula as in Proposition~\ref{prop:exp-lyapunov}. The approximations increase and remain below a Gaussian upper bound.
This use of monotonicity and Gaussian bounds is related to Kaniel and
Shinbrot \cite{KanielShinbrot1978}. We then apply Banach's theorem to the
map from the partner distribution to this solution
(Lemma~\ref{thm:boltzmann-local}).
The source estimates in Section~\ref{sec:bath-source} control the nonlinear term. Together with the linear contraction estimate \eqref{eq:linear-intersection-contraction}, they allow us to continue the solution for all time when \(\alpha\) is small.

\subsection{Mild Solutions on a Finite Interval}

On a time interval \(I\), a mild solution is a nonnegative
\(f\in C(I;L^1(W))\cap L^\infty_{\rm loc}(I;L^\infty(W))\) satisfying
\eqref{eq:boltzmann-mild} for every \(t_0\le t\) in \(I\).
\begin{lemma}[Collision estimates]\label{lem:gain-bounds}
For signed \(f,g\) for which the right-hand sides are finite,
\[
 \norm{\langle v\rangle^{-\gamma}Q(f,g)}_{L^\infty(W)}
 \le C\norm{f}_{L^\infty(W)}\norm{g}_{L^\infty(W)}
\]
and
\begin{align}
 \|Q(f,g)\|_{L^1(W)}
 &\le C\Bigl(
 \|\langle v\rangle^\gamma f\|_{L^1(W)}
 \|g\|_{L^\infty(W)}
 \notag\\
 &\hspace{35mm}
 +\|f\|_{L^\infty(W)}
 \|\langle v\rangle^\gamma g\|_{L^1(W)}
 \Bigr).
 \label{eq:collision-integrated-one}
\end{align}
Both bounds also hold for \(Q^+\) and \(Q^-\).

For almost every \(x\),
\begin{equation}\label{eq:collision-product-moment}
 \int_{\R}|Q(f,g)(x,v)|\dd v
 \le C
 \left(\int_{\R}\langle v\rangle^\gamma|f(x,v)|\dd v\right)
 \left(\int_{\R}\langle u\rangle^\gamma|g(x,u)|\dd u\right),
\end{equation}
and this also holds for \(Q^+\) and \(Q^-\). Finally,
\begin{equation}\label{eq:bath-collision-integrated-one}
 \|Q(f,M_{\tau(\cdot)})\|_{L^1(W)}
 \le C\|\langle v\rangle^\gamma f\|_{L^1(W)}.
\end{equation}
\end{lemma}

\begin{proof}
It suffices to take \(f,g\ge0\), since
\(\lvert Q^\pm(f,g)\rvert\le Q^\pm(\lvert f\rvert,\lvert g\rvert)\).

For the supremum estimate, use the gain formula
\eqref{eq:strong-gain-form}. Energy conservation and
\begin{equation}\label{eq:relative-speed-bound}
 |v-v_*|^\gamma\le C\langle v\rangle^\gamma\langle v_*\rangle^\gamma
\end{equation}
give
\[
\begin{aligned}
 W(v)Q^+(f,g)(x,v)
 &\le \|f\|_{L^\infty(W)}\|g\|_{L^\infty(W)}
 \int_{\R}|v-v_*|^\gamma e^{-a|v_*|^2}\dd v_*\\
 &\le C\langle v\rangle^\gamma
 \|f\|_{L^\infty(W)}\|g\|_{L^\infty(W)}.
\end{aligned}
\]
For the loss term, the same relative-speed bound gives
\begin{equation}\label{eq:collision-frequency-upper}
 |\nu[g](x,v)|
 \le C\langle v\rangle^\gamma
       \int_{\R}\langle u\rangle^\gamma|g(x,u)|\dd u
  \le C\langle v\rangle^\gamma\|g\|_{L^\infty(W)}.
\end{equation}
Multiplying by \(|f|\) proves the supremum bound for the
loss. Multiplication by \(W|f|\) and integration also gives its
\(L^1(W)\) bound.

For the \(L^1(W)\) gain estimate, fix \(x\).
Energy conservation gives
\[
 W(v')=W(v)W(v_*)e^{-a|v_*'|^2}.
\]
The weak form \eqref{eq:weak-gain-form} gives
\[
\begin{aligned}
 &\int_{\R}W(w)Q^+(f,g)(x,w)\dd w\\
 &\quad=\int_{\R}\int_{\R}\int_{\mathbb S^2}
 b_n|v-v_*|^\gamma W(v)f(x,v)W(v_*)g(x,v_*)
 e^{-a|v_*'|^2}\dd n\dd v_*\dd v .
\end{aligned}
\]
Changing from \(n\) to \(\sigma\) as in
\eqref{eq:weak-gain-form}, we have
\begin{equation}\label{eq:angular-gaussian-symmetry}
 \int_{\mathbb S^2}b_n e^{-a|v_*'|^2}\,\dd n
 =\frac1{4\pi}\int_{\mathbb S^2}
  e^{-\frac a4|v+v_*-|v-v_*|\sigma|^2}\,\dd\sigma.
\end{equation}
This angular average is symmetric in \(v,v_*\).
Split the velocity integral into \(|v|\le|v_*|\) and its complement.
On the first region, use
\[
 W(v)f(x,v)\le\|f\|_{L^\infty(W)},
 \qquad
 |v-v_*|^\gamma\le C\langle v_*\rangle^\gamma.
\]
The remaining Gaussian integral is bounded uniformly in
\(v_*\). Put \(R=|v_*|\) and let \(\theta\) be the angle
between \(n\) and \(v_*\). Then
\(|v_*'|^2=R^2\sin^2\theta+(v\cdot n)^2\).
Using \(b_n\le1/(2\pi)\), integration over \(|v|\le R\)
gives at most \(CR^2e^{-aR^2\sin^2\theta}\), and
\[
 R^2\int_0^\pi e^{-aR^2\sin^2\theta}\sin\theta\,\dd\theta\le2/a.
\]
The contribution from \(|v|\le|v_*|\) is therefore at most
\[
 C\|f\|_{L^\infty(W)}
 \int_{\R}\langle v_*\rangle^\gamma
 W(v_*)g(x,v_*)\,\dd v_*.
\]
The other region gives the same bound with \(f\) and \(g\)
interchanged. Adding and integrating in \(x\) proves
\eqref{eq:collision-integrated-one} for the gain.
The triangle inequality gives both bounds for \(Q=Q^+-Q^-\).

Testing \eqref{eq:weak-gain-form} with \(1\), then using
\eqref{eq:relative-speed-bound}, proves \eqref{eq:collision-product-moment}
for both gain and loss. For \eqref{eq:bath-collision-integrated-one},
the gain is bounded by \eqref{eq:bath-weighted-gain-moment}, and the
loss by \eqref{eq:hard-frequency-comparison}.
\end{proof}

\paragraph{Mild and distributional solutions.}
\label{par:mild-distributional}
On a bounded time interval \(I\), a nonnegative
\(f\in C(I;L^1(W))\cap L^\infty(I;L^\infty(W))\) is a distributional
solution of \eqref{eq:boltzmann} if and only if it satisfies
\eqref{eq:boltzmann-mild}.

Write the equation as \((\partial_t+v\cdot\nabla_x)f=G-af\), where \(G\) is the total gain and \(a\) the total collision frequency. The Gaussian upper bound gives uniformly bounded velocity moments of order \(\gamma\). Thus the frequency bounds and \eqref{eq:collision-product-moment} make both gain and loss integrable on \(I\times\T\times\R\), and give \(a(t,x,v)\le C\langle v\rangle^\gamma\).

For a distributional solution, set \(F(t,y,v)=f(t,y+tv,v)\). Changing variables in the weak equation removes transport and gives \(\partial_tF+\widetilde aF=\widetilde G\), where tildes mean evaluation at \((t,y+tv,v)\). Since \(\widetilde G-\widetilde aF\) is integrable, \(F\) has an absolutely continuous representative in time for almost every \((y,v)\). The bound on \(a\) makes the integrating factor well defined, and the product rule gives
\[
 \frac{\dd}{\dd t}\left(e^{\int_{t_0}^t\widetilde a(r)\dd r}F(t)\right)
 =e^{\int_{t_0}^t\widetilde a(r)\dd r}\widetilde G(t).
\]
Integrating and dividing by the exponential gives \eqref{eq:boltzmann-mild} after returning to \(x=y+tv\).

Conversely, the mild formula is absolutely continuous along almost every such path. Differentiating gives \(\partial_tF+\widetilde aF=\widetilde G\). Integration against a smooth test function and the change back to \(x\) recover the distributional equation. Continuity in \(L^1\) gives the mild identity at every pair of times \(t_0\le t\), as an equality in \(L^1\).

\begin{lemma}[Construction on a finite interval]
\label{thm:boltzmann-local}
Fix \(A,S<\infty\) and \(H>0\).  There are
\(0<\bar\alpha\le1/2\) and \(C<\infty\), depending only on
\(A,S,H\) and the model parameters, such that, if
\(0\le\alpha\le\bar\alpha\) and \(f_{\rm in}\) is a probability density satisfying
\[
 \norm{f_{\rm in}}_{L^1(W)}\le A,
 \qquad \norm{f_{\rm in}}_{L^\infty(W)}\le S,
\]
then there is a nonnegative mild solution on \([0,H]\) with
\[
 \sup_{t\le H}\bigl(\|f(t)\|_{L^1(W)}+\|f(t)\|_{L^\infty(W)}\bigr)
 +\int_0^H\|\langle v\rangle^\gamma f(t)\|_{L^1(W)}\dd t\le C.
\]
It is the unique nonnegative mild solution with this datum.
\end{lemma}

\begin{proof}
For fixed \(g\ge0\), consider the linear equation
\[
 \partial_t f_g+v\cdot\nabla_xf_g
 =(1-\alpha)Q(f_g,M_{\tau(x)})+\alpha Q(f_g,g),
 \qquad f_g(0)=f_{\rm in}.
\]
Choose \(c,C_0\) as in Lemma~\ref{lem:bath-drift} and \(C\) as in
Lemma~\ref{lem:gain-bounds}. Put \(B=S e^{C_0H}\) and suppose that
\(\|g(t)\|_{L^\infty(W)}\le B\) for almost every \(t\le H\).
If \(\alpha\le1/2\) and \(C\alpha B\le c/8\), then
{
\[
\begin{aligned}
 &(1-\alpha)Q(W^{-1},M_{\tau(x)})+\alpha Q(W^{-1},g)\\
 &\quad\le
 \bigl[C_0-\{(1-\alpha)c-C\alpha B\}\langle v\rangle^\gamma\bigr]W^{-1}
 \le C_0W^{-1}.
\end{aligned}
\]
}
Thus \(S e^{C_0t}W^{-1}\) is a supersolution. Comparison gives
\[
 0\le f_g(t,x,v)\le S e^{C_0t}W^{-1}(v),
 \qquad \sup_{t\le H}\|f_g(t)\|_{L^\infty(W)}\le B.
\]

The Gaussian comparison only gives \(Wf_g\le B\), which need not
be integrable in velocity. Formally, mass conservation and testing
with \(W\) give, by Lemmas~\ref{lem:bath-drift} and~\ref{lem:gain-bounds},
{
\begin{equation}\label{eq:frozen-weighted-moment}
 \|f_g(t)\|_{L^1(W)}
 +\{(1-\alpha)c-C\alpha B\}
   \int_0^t\|\langle v\rangle^\gamma f_g(s)\|_{L^1(W)}\dd s
 \le A+C_0t+C\alpha B
   \int_0^t\|\langle v\rangle^\gamma g(s)\|_{L^1(W)}\dd s.
\end{equation}
}
Estimating the nonlinear gain requires the stronger moments of both
\(f_g\) and \(g\). The \(f_g\) term is absorbed by the bath loss
in \eqref{eq:frozen-weighted-moment}, leaving a coefficient of
at least \(c/4\). For \(g\), we choose a bound that will also
hold for \(f_g\). Set
\[
 D=\frac{8(A+C_0H)}c.
\]
The condition \(C\alpha B\le c/8\) gives
\(C\alpha BD\le A+C_0H\). If
\(\int_0^H\|\langle v\rangle^\gamma g(t)\|_{L^1(W)}\dd t\le D\),
the right side of
\eqref{eq:frozen-weighted-moment} is therefore at most
\(2(A+C_0H)\), and
\[
 \sup_{t\le H}\|f_g(t)\|_{L^1(W)}\le2(A+C_0H),
 \qquad
 \int_0^H\|\langle v\rangle^\gamma f_g(t)\|_{L^1(W)}\dd t\le D.
\]

Work in the space
\[
 Y=L^\infty((0,H)\times\T\times\R;W),
 \qquad \|g\|_Y=\operatorname*{ess\,sup}_{0<t<H}
                     \|g(t)\|_{L^\infty(W)},
\]
with the subset
\[
 \mathcal K=\left\{g\in Y:g\ge0,\ \|g\|_Y\le B,\quad
 \int_0^H\|\langle v\rangle^\gamma g(t)\|_{L^1(W)}\dd t\le D\right\}.
\]
The set \(\mathcal K\) is nonempty and, by Fatou's lemma, closed
in \(Y\). 

For fixed \(g\in\mathcal K\), set \(f^{(0)}=0\) and let
\(f^{(n+1)}\) solve, in the mild sense,
{
\begin{equation}\label{eq:frozen-gain-iteration}
\begin{aligned}
 &(\partial_t+v\cdot\nabla_x)f^{(n+1)}
 +\bigl((1-\alpha)\nu[M_{\tau(x)}]+\alpha\nu[g]\bigr)f^{(n+1)}\\
 &\qquad=(1-\alpha)Q^+(f^{(n)},M_{\tau(x)})
              +\alpha Q^+(f^{(n)},g),\\
 &f^{(n+1)}(0)=f_{\rm in}.
\end{aligned}
\end{equation}
}
As in the linear construction in Section~\ref{sec:linear-construction}, the loss coefficient is independent of \(n\), so \(0\le f^{(n)}\le f^{(n+1)}\). Comparison with \(S e^{C_0t}W^{-1}\) bounds the iterates uniformly in \(n\).
We can pass to the limit in the characteristic formula by monotone convergence,
obtaining a solution \(f_g\le S e^{C_0t}W^{-1}\).
Estimate~\eqref{eq:collision-product-moment} permits
integration of the equation in \(x,v\), which proves
\(\int_{\T\times\R}f_g(t)=1\); continuity in \(L^1\)
also follows from the characteristic formula.

To justify \eqref{eq:frozen-weighted-moment} without assuming
the moment is finite, test with \(W_N=\min\{W,N\}\), as in
Proposition~\ref{prop:exp-lyapunov}.
Since \(W_N(v')-W_N(v)\le0\) on \(\{W\ge N\}\),
{
\[
\begin{aligned}
 \int_{\T\times\R}W_NQ(f_g,g)
 &\le \int_{\T\times\R}WQ^+(f_g{\bf1}_{\{W<N\}},g)\\
 &\le CB\left(
 \int_{\T\times\R}\langle v\rangle^\gamma
 W{\bf1}_{\{W<N\}}f_g
 +\|\langle v\rangle^\gamma g\|_{L^1(W)}
 \right).
\end{aligned}
\]
}
For the bath term, apply the calculation in
Proposition~\ref{prop:exp-lyapunov}. Adding the last inequality
leaves the coefficient \((1-\alpha)c-C\alpha B\ge c/4\) on the
restricted moment integral. We can therefore integrate in time and
send \(N\to\infty\). By monotone convergence, this proves
\eqref{eq:frozen-weighted-moment} and the moment bounds.

The moment bound and
\eqref{eq:collision-integrated-one}--\eqref{eq:bath-collision-integrated-one}
show that the source
\(\alpha\{Q(f_g,g)-Q(f_g,M_{\tau(\cdot)})\}\) belongs to
\(L^1((0,H);L^1(W))\). Lemma~\ref{lem:bath-source-duhamel}
applied to \(f_g-P_t^\tau f_{\rm in}\) therefore gives the
representation
\begin{equation}\label{eq:frozen-bath-duhamel}
 f_g(t)=P_t^\tau f_{\rm in}
 +\alpha\int_0^tP_{t-s}^\tau
 \{Q(f_g(s),g(s))-Q(f_g(s),M_{\tau(\cdot)})\}\,\dd s.
\end{equation}
Moreover, \(f_g\in C([0,H];L^1(W))\) by
Lemma~\ref{lem:bath-source-duhamel}. Together with the bounds above,
this proves that \(g\mapsto f_g\) maps \(\mathcal K\) into itself.

For \(g,h\in\mathcal K\), put \(u=f_g-f_h\).
The initial terms in \eqref{eq:frozen-bath-duhamel} cancel.
By bilinearity, the source for \(u\) is
\[
 \alpha\{Q(u,g)+Q(f_h,g-h)-Q(u,M_{\tau(\cdot)})\}.
\]
Both \(g\) and \(f_h\) are bounded by \(B\) in
\(L^\infty(W)\). The collision estimate and
\eqref{eq:bath-source-duhamel} therefore give
\[
 \|u\|_Y\le C\alpha\bigl(\|u\|_Y+\|g-h\|_Y\bigr).
\]
Here \(C\) depends only on \(A,S,H\) and the model parameters.
Choose \(C\alpha\le1/3\) and move the term containing \(u\)
to the left:
\[
 \|f_g-f_h\|_Y
 \le\frac{C\alpha}{1-C\alpha}\|g-h\|_Y
 \le\frac12\|g-h\|_Y.
\]
By Banach's theorem there is a unique \(g\in\mathcal K\) with
\(g=f_g\). This proves existence of a mild solution of \eqref{eq:boltzmann}.

Let \(f,g\) be two nonnegative mild solutions with the same
datum on \([0,H]\). The definition of a mild solution gives an
\(L^\infty(W)\) bound \(B_1\) for both solutions on this interval. Put \(w=f-g\). Subtracting the equations and using the symmetry
\(Q^+(g,w)=Q^+(w,g)\) established from
\eqref{eq:weak-gain-form} gives
\[
\begin{aligned}
 &(\partial_t+v\cdot\nabla_x)w+
 \bigl((1-\alpha)\nu[M_{\tau(x)}]+\alpha\nu[f]\bigr)w\\
 &\qquad=(1-\alpha)Q^+(w,M_{\tau(x)})
 +\alpha Q^+(w,f+g)-\alpha g\nu[w].
\end{aligned}
\]
We use \(W^{1/2}\) because
\(W^{1/2}(f+g)\le2B_1e^{-a|v|^2/2}\) remains integrable after
multiplication by \(\langle v\rangle^\gamma\). Thus the collision terms
are integrable with this weight without first using the stronger moment
bound. Integrate the inequality for \(|w|\) against \(W^{1/2}\), using
the angular calculation in \eqref{eq:bath-weighted-gain-moment} for
the gain.
The factor \(e^{a|v_*|^2/2}\) in that calculation is integrable
against both \(M_{\tau(x)}\) and \(f+g\le2B_1W^{-1}\).
Together with the bath loss and \(1-\alpha\ge1/2\), this gives,
for almost every \(t\),
{
\[
\begin{aligned}
 \frac{\dd}{\dd t}\|w(t)\|_{L^1(W^{1/2})}
 &\le\int_{\T\times\R}
 \left(\frac{C_{B_1}}{1+|v|}-\frac{c_{\rm bath}}2\right)
 \langle v\rangle^\gamma W^{1/2}|w|\,\dd v\dd x\\
 &\quad+C_{B_1}\int_{\T\times\R}
 \langle v\rangle^\gamma|w|\,\dd v\dd x\\
 &\le C_{B_1}\|w(t)\|_{L^1(W^{1/2})}.
\end{aligned}
\]
}
Here the contribution of \(g\nu[w]\) is estimated using
\(W^{1/2}g\le B_1e^{-a|v|^2/2}\).
In the first integral, the coefficient is negative for large
\(|v|\) and bounded on the remaining ball.
The second integral is bounded by
\(C_{B_1}\|w\|_{L^1(W^{1/2})}\), since
\(\langle v\rangle^\gamma\le CW^{1/2}\).
As \(w(0)=0\), Gr\"onwall's inequality proves uniqueness.
\end{proof}

\subsection{Global Continuation}

Choose \(t_c\) so that \eqref{eq:linear-intersection-contraction}
holds with \(\theta=1/4\). To repeat the local construction on intervals of length \(t_c\) with
\(\alpha\) fixed, we need the same norm bound at the end of each interval. We first estimate the linear part by applying
\eqref{eq:linear-intersection-contraction} to the difference of the
probability densities \(f\) and \(M_{\tau_-}\). Together with
\eqref{eq:linear-L1-uniform} and \eqref{eq:linear-Linf-uniform} for
\(P_{t_c}^\tau M_{\tau_-}\), this gives \(K>0\), uniform over baths
satisfying \eqref{eq:tau-bounds}, such that
\begin{equation}\label{eq:affine-bath-return}
\begin{aligned}
 \|P_{t_c}^\tau f\|_{L^1(W)}+\|P_{t_c}^\tau f\|_{L^\infty(W)}
 \le\frac14\bigl(\|f\|_{L^1(W)}+\|f\|_{L^\infty(W)}\bigr)+K
\end{aligned}
\end{equation}
for every probability density \(f\in L^1(W)\cap L^\infty(W)\).
Increase \(K\), if needed, so that
\(\|M_{\tau_-}\|_{L^1(W)}+\|M_{\tau_-}\|_{L^\infty(W)}\le K\). This ensures that \(\mathcal R\) contains \(M_{\tau_-}\), and is therefore nonempty, for the radii chosen below.

\begin{lemma}[Continuation to \(t_c\)]
\label{lem:interval-continuation}
There is \(R_*>0\) such that, for every \(R\ge R_*\), there are
\(\bar\alpha>0\) and \(C<\infty\) with the following property.  If
\(0\le\alpha\le\bar\alpha\) and \(f_{\rm in}\) is a probability density satisfying
\[
 \|f_{\rm in}\|_{L^1(W)}+\|f_{\rm in}\|_{L^\infty(W)}\le R,
\]
then the mild solution exists on \([0,t_c]\),
\begin{equation}\label{eq:interval-bound}
 \sup_{0\le t\le t_c}
 \bigl(\|f(t)\|_{L^1(W)}+\|f(t)\|_{L^\infty(W)}\bigr)\le C,
\end{equation}
and
\begin{equation}\label{eq:interval-return}
 \|f(t_c)\|_{L^1(W)}+\|f(t_c)\|_{L^\infty(W)}\le\frac{3R}{4}.
\end{equation}
The solution is unique among nonnegative mild solutions with this datum.
\end{lemma}
\begin{proof}
Set \(R_*=4K\) and fix \(R\ge R_*\). Apply
Lemma~\ref{thm:boltzmann-local} with \(A=S=R\) and \(H=t_c\).
For sufficiently small \(\alpha\), it gives a unique solution on \([0,t_c]\)
and a bound \(C_R\), independent of \(\alpha\), such that
\begin{equation}\label{eq:continuation-integrated-bound}
 \sup_{t\le t_c}\bigl(\|f(t)\|_{L^1(W)}+\|f(t)\|_{L^\infty(W)}\bigr)
 +\int_0^{t_c}\|\langle v\rangle^\gamma f(t)\|_{L^1(W)}\dd t\le C_R.
\end{equation}
Duhamel's formula around the linear evolution gives
\[
 f(t)=P_t^\tau f_{\rm in}
 +\alpha\int_0^tP_{t-s}^\tau
 \{Q(f(s),f(s))-Q(f(s),M_{\tau(\cdot)})\}\dd s .
\]
Lemma~\ref{lem:gain-bounds} and
\eqref{eq:continuation-integrated-bound} give
\[
\begin{aligned}
 &\operatorname*{ess\,sup}_{s\le t_c}
 \|\langle v\rangle^{-\gamma}
   \{Q(f,f)-Q(f,M_{\tau(\cdot)})\}(s)\|_{L^\infty(W)}\le C_R,\\
 &\int_0^{t_c}
 \|\{Q(f,f)-Q(f,M_{\tau(\cdot)})\}(s)\|_{L^1(W)}\dd s\le C_R.
\end{aligned}
\]
The first bound uses \(\sup_s\|f(s)\|_{L^\infty(W)}\);
the second also uses the moment integral in
\eqref{eq:continuation-integrated-bound}.
Lemma~\ref{lem:bath-source-duhamel} therefore gives
\[
 \sup_{t\le t_c}
 \bigl(\|f(t)-P_t^\tau f_{\rm in}\|_{L^1(W)}
 +\|f(t)-P_t^\tau f_{\rm in}\|_{L^\infty(W)}\bigr)
 \le C_R\alpha .
\]
At \(t_c\), \eqref{eq:affine-bath-return} gives
\[
 \|f(t_c)\|_{L^1(W)}+\|f(t_c)\|_{L^\infty(W)}
 \le R/4+K+C_R\alpha\le3R/4,
\]
after decreasing \(\alpha\) so that \(C_R\alpha\le R/4\).
The bound on \([0,t_c]\) follows from
\eqref{eq:continuation-integrated-bound}.
Uniqueness follows from Lemma~\ref{thm:boltzmann-local}.
\end{proof}

\subsection{Contraction and the Stationary State}
\label{sec:nonlinear-contraction}

Since two probability densities have the same mass,
\eqref{eq:linear-intersection-contraction} contracts their difference.
We want to keep this contraction after adding self collisions, so we
estimate the nonlinear integral in terms of the initial difference.

\begin{lemma}[Difference up to \(t_c\)]
\label{prop:nonlinear-time-contraction}
Let \(f\) and \(g\) be mild solutions of unit mass on
\([0,t_c]\) satisfying
\[
 \sup_{0\le t\le t_c}
 \bigl(\|f(t)\|_{L^1(W)}+\|f(t)\|_{L^\infty(W)}
 +\|g(t)\|_{L^1(W)}+\|g(t)\|_{L^\infty(W)}\bigr)\le B.
\]
There is \(C<\infty\), depending on \(B\), such that, for sufficiently
small \(\alpha\),
\begin{align}
 &\begin{aligned}
 &\sup_{0\le t\le t_c}
 \bigl(\|f(t)-g(t)\|_{L^1(W)}+\|f(t)-g(t)\|_{L^\infty(W)}\bigr)\\
 &\qquad\le C\bigl(\|f(0)-g(0)\|_{L^1(W)}
 +\|f(0)-g(0)\|_{L^\infty(W)}\bigr),
 \end{aligned}
 \label{eq:nonlinear-time-stability}\\
 &\begin{aligned}
 &\|f(t_c)-g(t_c)\|_{L^1(W)}+\|f(t_c)-g(t_c)\|_{L^\infty(W)}\\
 &\qquad\le\frac12\bigl(\|f(0)-g(0)\|_{L^1(W)}
 +\|f(0)-g(0)\|_{L^\infty(W)}\bigr).
 \end{aligned}
 \label{eq:nonlinear-time-contraction}
\end{align}
\end{lemma}

\begin{proof}
Set \(w=f-g\).

Lemma~\ref{thm:boltzmann-local}, with \(A=S=B\) and \(H=t_c\),
applies to both initial data. By uniqueness, its moment bounds hold
for \(f\) and \(g\). Thus, for sufficiently small \(\alpha\),
\begin{equation}\label{eq:solution-integrated-bound}
 \int_0^{t_c}
 \|\langle v\rangle^\gamma f(t)\|_{L^1(W)}\dd t
 +\int_0^{t_c}
 \|\langle v\rangle^\gamma g(t)\|_{L^1(W)}\dd t
 \le C_B.
\end{equation}
Since \(|w|\le f+g\), Lemma~\ref{lem:gain-bounds} makes
the collision terms below integrable in \(L^1(W)\) over time.
Bilinearity and Lemma~\ref{lem:bath-source-duhamel} give
\[
 w(t)=P_t^\tau w(0)+\alpha\int_0^tP_{t-s}^\tau
 \bigl[Q(w,f)+Q(g,w)-Q(w,M_{\tau(\cdot)})\bigr](s)\,\dd s.
\]

The Gaussian bounds on \(f\) and \(g\) give
\begin{equation}\label{eq:mixed-source-infinity}
 \left\|\langle v\rangle^{-\gamma}
 \{Q(w,f)+Q(g,w)-Q(w,M_{\tau(\cdot)})\}\right\|_{L^\infty(W)}
 \le C_B\|w\|_{L^\infty(W)}.
\end{equation}
Apply Lemma~\ref{lem:bath-Linf-boundedness} to the initial
term and \eqref{eq:bath-source-duhamel} to the integral:
\[
 \sup_{0\le t\le t_c}
 \|w(t)\|_{L^\infty(W)}
 \le C\|w(0)\|_{L^\infty(W)}
 +C_B\alpha\sup_{0\le t\le t_c}
 \|w(t)\|_{L^\infty(W)}.
\]
For \(C_B\alpha\le1/2\), moving the last term to the left gives
\begin{equation}\label{eq:difference-infinity-bound}
 \sup_{0\le t\le t_c}
 \|w(t)\|_{L^\infty(W)}
 \le C_B\|w(0)\|_{L^\infty(W)}.
\end{equation}

For the weighted \(L^1\) norm, the mixed collision estimate gives
\begin{align}
 &\|Q(w,f)+Q(g,w)-Q(w,M_{\tau(\cdot)})\|_{L^1(W)}
 \notag\\
 &\quad\le C_B\|\langle v\rangle^\gamma w\|_{L^1(W)}
 +C\|w\|_{L^\infty(W)}
 \Bigl(\|\langle v\rangle^\gamma f\|_{L^1(W)}
 +\|\langle v\rangle^\gamma g\|_{L^1(W)}\Bigr).
 \label{eq:mixed-source-one}
\end{align}
By \eqref{eq:solution-integrated-bound} and
\eqref{eq:difference-infinity-bound}, the time integral of the
second term is at most \(C_B\|w(0)\|_{L^\infty(W)}\).
The bound \(|w|\le f+g\) proves integrability of the first term,
but loses the dependence on \(w(0)\). We keep the stronger moment
of \(w\) on the left in the linear source estimate so that we can
bound its integral in terms of \(w(0)\). Applying \eqref{eq:bath-integrated-moment} to
\(P_t^\tau w(0)\) and
\eqref{eq:bath-source-one-integrated} to the integral term gives
\begin{align*}
 &\sup_{0\le t\le t_c}\|w(t)\|_{L^1(W)}
 +\int_0^{t_c}\|\langle v\rangle^\gamma w(t)\|_{L^1(W)}\dd t\\
 &\quad\le C\|w(0)\|_{L^1(W)}
 +C_B\alpha\int_0^{t_c}
 \|\langle v\rangle^\gamma w(t)\|_{L^1(W)}\dd t\\
 &\qquad+C_B\alpha\|w(0)\|_{L^\infty(W)}.
\end{align*}
Decrease \(\alpha\), if necessary, so that \(C_B\alpha\le1/2\).
Moving the integral on the right to the left gives
\begin{align}
 &\sup_{0\le t\le t_c}\|w(t)\|_{L^1(W)}
 +\int_0^{t_c}
 \|\langle v\rangle^\gamma w(t)\|_{L^1(W)}\dd t
 \notag\\
 &\qquad\le C_B\left(
 \|w(0)\|_{L^1(W)}+\|w(0)\|_{L^\infty(W)}\right).
 \label{eq:difference-integrated-bound}
\end{align}
Together with \eqref{eq:difference-infinity-bound}, this proves
\eqref{eq:nonlinear-time-stability}.

Apply Lemma~\ref{lem:bath-source-duhamel} using
\eqref{eq:mixed-source-infinity}--\eqref{eq:difference-integrated-bound}
to obtain
\begin{align*}
 &\|w(t_c)-P_{t_c}^\tau w(0)\|_{L^1(W)}
 +\|w(t_c)-P_{t_c}^\tau w(0)\|_{L^\infty(W)}\\
 &\qquad\le C_B\alpha\bigl(\|w(0)\|_{L^1(W)}+\|w(0)\|_{L^\infty(W)}\bigr).
\end{align*}
Since \(\int w(0)=0\), our choice of \(t_c\) lets us apply
\eqref{eq:linear-intersection-contraction} with \(\theta=1/4\), giving
\begin{align*}
 \|w(t_c)\|_{L^1(W)}+\|w(t_c)\|_{L^\infty(W)}
 \le\left(\frac14+C_B\alpha\right)
 \bigl(\|w(0)\|_{L^1(W)}+\|w(0)\|_{L^\infty(W)}\bigr).
\end{align*}
Taking \(C_B\alpha\le1/4\) proves
\eqref{eq:nonlinear-time-contraction}.
\end{proof}
\begin{proof}[Proof of Theorem~\ref{thm:intro-relaxation}]
Fix \(t_c\) as in \eqref{eq:affine-bath-return} and take
\(R\ge R_*\). Let \(C_{\rm cont}(R)\) be the bound in
Lemma~\ref{lem:interval-continuation} for this radius.
Choose \(\alpha_0>0\) small enough that
Lemma~\ref{lem:interval-continuation} and
Lemma~\ref{prop:nonlinear-time-contraction} with
\(B=2C_{\rm cont}(R)\) both apply.
By \eqref{eq:interval-return}, \(f(t_c)\in\mathcal R\).
Inductively, \(f(nt_c)\in\mathcal R\) for every integer \(n\ge0\),
and \eqref{eq:interval-bound} holds on each interval
\([nt_c,(n+1)t_c]\). This constructs the solution for all \(t\ge0\).
The map \(f_{\rm in}\mapsto f(t_c)\) preserves \(\mathcal R\) and,
by \eqref{eq:nonlinear-time-contraction}, has Lipschitz constant
\(1/2\). Since \(\mathcal R\) is closed in
\(L^1(W)\cap L^\infty(W)\) with the sum norm, Banach's theorem
provides a unique fixed point \(F_\alpha\).

For this fixed point, \(f(t_c)=f(0)=F_\alpha\).
Uniqueness implies \(f(t+t_c)=f(t)\) for every \(t\ge0\).
Fix \(s\ge0\) and apply \eqref{eq:nonlinear-time-contraction} to
\(f(t+s)\) and \(f(t)\) on \([0,t_c]\).
By periodicity their distance is the same at both endpoints, so
contraction forces this distance to be zero.
We conclude that \(f(s)=F_\alpha\) for all \(s\ge0\).

For \(f_{\rm in}\in\mathcal R\), iteration of
\eqref{eq:nonlinear-time-contraction} gives, for every integer \(n\ge0\),
\[
\begin{aligned}
 &\|f(nt_c)-F_\alpha\|_{L^1(W)}
 +\|f(nt_c)-F_\alpha\|_{L^\infty(W)}\\
 &\qquad\le2^{-n}\bigl(\|f_{\rm in}-F_\alpha\|_{L^1(W)}
 +\|f_{\rm in}-F_\alpha\|_{L^\infty(W)}\bigr).
\end{aligned}
\]
For \(nt_c\le t<(n+1)t_c\), apply
\eqref{eq:nonlinear-time-stability} from time \(nt_c\).
The last display then proves the decay in both norms.
For the moment integral, apply \eqref{eq:difference-integrated-bound}
on each later interval and sum. This gives the bound on the integral from \(t\) to infinity in Theorem~\ref{thm:intro-relaxation}.

By the \hyperref[par:mild-distributional]{characteristic argument},
a stationary distributional state in \(\mathcal R\) is also a mild
solution. It is therefore a fixed point of the time-\(t_c\) map and
equals \(F_\alpha\).
The stationary state for the smaller class also belongs to any
larger class. If \(\alpha\) satisfies the smallness conditions for
both radii, uniqueness in the larger class proves that the two
stationary states coincide.
\end{proof}
\section{Properties of Stationary States}
\label{sec:balance-laws}

Collisions with the reservoir change the gas momentum and
energy, whereas self collisions preserve both. For Maxwell
molecules \((\gamma=0)\), we can compute these changes from the
density, momentum and energy of \(f\).

\subsection{Balance Laws}

Write
\[
 m_f=\int_{\R}vf\,\dd v,\qquad
 P_f=\int_{\R}v\otimes v\,f\,\dd v,\qquad
 q_f=\frac12\int_{\R}v|v|^2f\,\dd v,
\]%
and \((\nabla_x\cdot P_f)_i=\partial_j(P_f)_{ij}\), with repeated
indices summed. The following balance laws apply in particular to the
stationary state of Theorem~\ref{thm:boltzmann-contraction}.

\begin{proposition}[Balance laws]\label{prop:balance}
Let \(\gamma=0\), \(0\leq\alpha<1\), and let \(f\) be a nonnegative
stationary distributional solution of \eqref{eq:boltzmann} with
\(f\leq Ce^{-c|v|^2}\) for some \(c,C>0\). Then, in the sense of
distributions,
\begin{equation}\label{eq:balance-laws}
 \nabla_x\cdot m_f=0,\qquad
 \nabla_x\cdot P_f=-\frac{1-\alpha}2m_f,\qquad
 \nabla_x\cdot q_f=\frac{1-\alpha}4\bigl(3\tau\rho_f-e_f\bigr).
\end{equation}%
Consequently \(\partial_i\partial_j(P_f)_{ij}=0\),
\(\int_\T m_f\,\dd x=0\) and \(\int_\T e_f\,\dd x=3\int_\T\tau\rho_f\,\dd x\).
If \(\rho_f>0\) almost everywhere, the last identity reads
\begin{equation}\label{eq:stationary-temperature-deficit}
 \int_\T\rho_f(\tau-T_f)\,\dd x=\int_\T\frac{|m_f|^2}{3\rho_f}\,\dd x\geq0.
\end{equation}%

\end{proposition}

\begin{proof}
The Gaussian bound allows us to test the stationary equation with
\(\varphi(x)\psi(v)\), where \(\varphi\in C^\infty(\T)\) and
\(\psi\in\{1,v_j,|v|^2\}\). This gives
\begin{equation}\label{eq:moment-balance}
  \nabla_x\cdot\int_{\R} v\psi f\,\dd v
  =(1-\alpha)\int_{\R}\psi Q(f,M_{\tau(x)})\,\dd v
   +\alpha\int_{\R}\psi Q(f,f)\,\dd v .
\end{equation}%
Self collisions conserve mass, momentum, and energy, so their
integrals in \eqref{eq:moment-balance} vanish for these tests.
For collisions with the reservoir, use the last integral in
\eqref{eq:weak-gain-form}. Since
\(\int_{\mathbb S^2}\sigma\,\dd\sigma=0\),
\begin{equation}\label{eq:angular-moment-increments}
\begin{aligned}
 \frac1{4\pi}\int_{\mathbb S^2}(v'-v)\dd\sigma
 &=\frac{v_*-v}{2},\\
 \frac1{4\pi}\int_{\mathbb S^2}(|v'|^2-|v|^2)\dd\sigma
 &=\frac{|v_*|^2-|v|^2}{2}.
\end{aligned}
\end{equation}
Integrating \eqref{eq:angular-moment-increments} against
\(f(v)M_{\tau(x)}(v_*)\), using the zero mean and energy
\(3\tau\) of the reservoir Maxwellian, gives
\begin{equation}\label{eq:maxwell-bath-moments}
 \int_{\R}vQ(f,M_{\tau(x)})\,\dd v=-\frac12m_f,
 \qquad
 \int_{\R}|v|^2Q(f,M_{\tau(x)})\,\dd v=\frac12(3\tau\rho_f-e_f).
\end{equation}%
Inserting \eqref{eq:maxwell-bath-moments} in
\eqref{eq:moment-balance} with \(\psi=1\), \(v_j\) and \(|v|^2/2\) proves
\eqref{eq:balance-laws}. Taking the divergence of the second identity and
using the first gives \(\partial_i\partial_j(P_f)_{ij}=0\); integrating the
second and third identities over the torus gives the other two
consequences; and \(e_f=3\rho_fT_f+|m_f|^2/\rho_f\) gives
\eqref{eq:stationary-temperature-deficit}.
\end{proof}

If \(\tau=\tau(x_1)\) and \(f=f(x_1,v)\),
\eqref{eq:balance-laws} reduces to equations in \(x_1\).

\begin{corollary}[Dependence on one coordinate]\label{cor:planar-balance}
In Proposition~\ref{prop:balance}, suppose in addition that
\(\tau=\tau(x_1)\) and \(f=f(x_1,v)\). Then
\((m_f)_1\equiv0\), \((P_f)_{11}\) is a positive constant, \(\rho_f>0\)
almost everywhere, and \((q_f)_1\) is Lipschitz with
\[
 \partial_{x_1}(q_f)_1=\frac{1-\alpha}4\bigl(3\tau\rho_f-e_f\bigr).
\]
If \(f\) is also invariant under rotations and reflections of
\((v_2,v_3)\), as \(F_\alpha\) is
(Section~\ref{sec:symmetry-reductions}), then \(m_f=0\) and
\begin{equation}\label{eq:planar-energy-exchange}
 \partial_{x_1}(q_f)_1=\frac{3(1-\alpha)}4\rho_f(\tau-T_f).
\end{equation}
\end{corollary}

\begin{proof}
For functions of \(x_1\), the first two identities in
\eqref{eq:balance-laws} read \(\partial_{x_1}(m_f)_1=0\) and
\(\partial_{x_1}(P_f)_{11}=-\frac{1-\alpha}2(m_f)_1\). Thus \((m_f)_1\) is
constant. Integrating the second identity over \(\mathbb T^1\) gives
\((m_f)_1=0\), and \((P_f)_{11}\) is constant. It is positive, since a probability
density cannot be supported on the plane \(\{v_1=0\}\).
Then \(\int_\R v_1^2f(x_1,v)\,\dd v>0\) gives
\(\rho_f>0\) almost everywhere. The Gaussian bound makes \(\rho_f\) and
\(e_f\) bounded, so \((q_f)_1\) is Lipschitz. Under the symmetries,
\((m_f)_2=(m_f)_3=0\), so \(e_f=3\rho_fT_f\).
\end{proof}

\paragraph{Energy exchange.}
For Maxwell molecules and \(f=F_\alpha\) with \(\tau=\tau(x_1)\),
\eqref{eq:planar-energy-exchange} shows that \((q_f)_1\) increases where
\(T_f<\tau\) and decreases where \(T_f>\tau\). Proposition~\ref{prop:stationary-entropy-balance} rules out a
constant current for a nonconstant reservoir.

\begin{figure}[!b]
  \centering
  \includegraphics[width=0.98\textwidth]
    {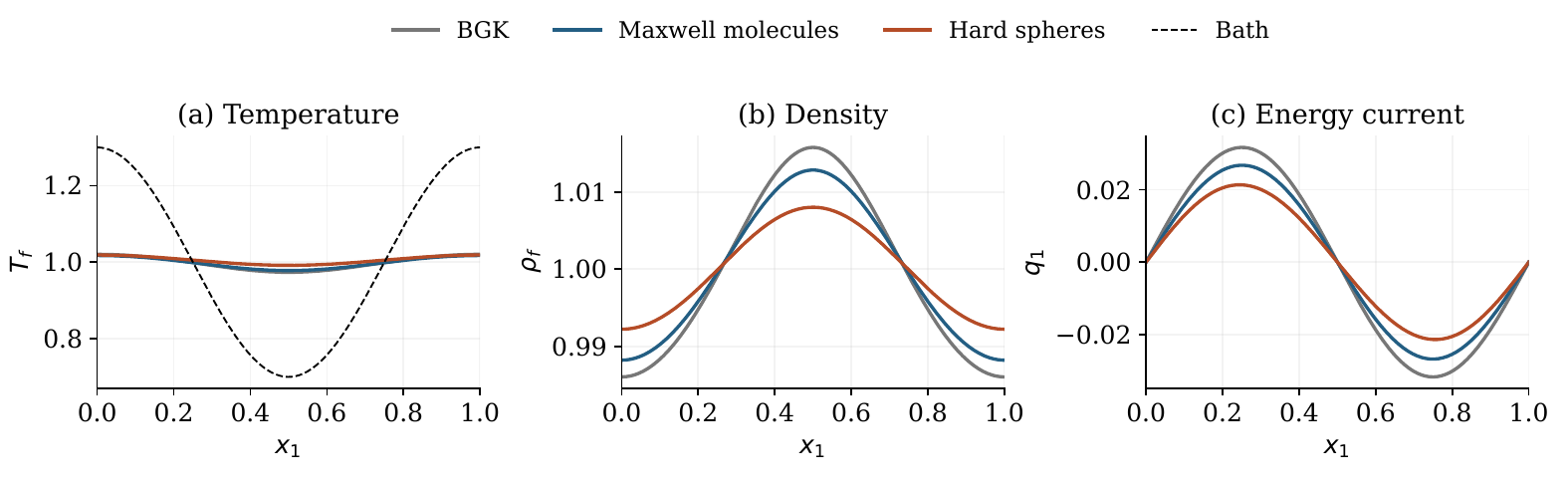}
  \caption{Computed stationary fields for BGK, Maxwell molecules,
  and hard spheres, with \(\tau(x_1)=1+0.3\cos(2\pi x_1)\)
  and \(\alpha=0.2\).
  (a) Gas and reservoir temperatures. (b) Density. (c) Energy current.
  Parameters and collision rates are given in
  Appendix~\ref{app:numerical-illustrations}.}
  \label{fig:physical-hard-sphere-fields}
\end{figure}

In Figure~\ref{fig:physical-hard-sphere-fields} we compare the
Boltzmann and BGK stationary fields for the same reservoir.

\subsection{The Energy Current as a Function of Velocity}
\label{sec:energy-current-velocity}

\begin{figure}[H]
  \centering
  \includegraphics[width=0.94\textwidth]
    {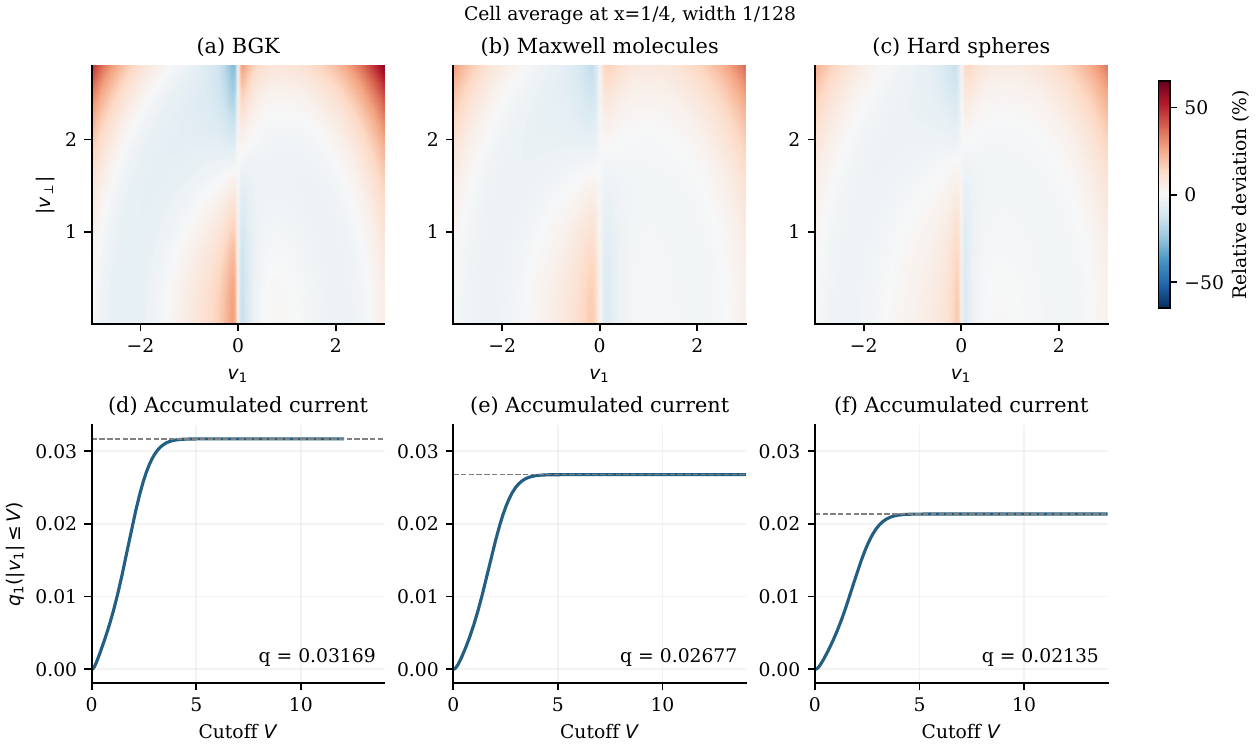}
  \caption{Velocity-space comparison near \(x_1=1/4\).
  Top: relative deviation from the local Maxwellian (\%), using the same
  colour scale in all three panels.
  Bottom: the cumulative current \(q_1(V)\) defined in
  \eqref{eq:cumulative-energy-current}.
  Dashed lines show the computed total currents.}
  \label{fig:stationary-velocity}
\end{figure}

Figure~\ref{fig:stationary-velocity} compares each computed distribution
with the local Maxwellian having the same density, mean velocity and temperature,
\[
 u_f=\rho_f^{-1}\int_{\R}vf\,\dd v,
 \qquad M_{\rm loc}=\rho_fM_{u_f,T_f}.
\]

For the displayed states, the transverse symmetries of
Section~\ref{sec:symmetry-reductions} give
\((m_f)_2=(m_f)_3=0\). Mass conservation gives
\(\partial_{x_1}(m_f)_1=0\), and the reflection symmetry
\((x_1,v_1)\mapsto(-x_1,-v_1)\) makes this constant odd.
Hence \(m_f=0\) and \(u_f=0\).

The local Maxwellian therefore carries no energy current, and
\[
 (q_f)_1=\frac12\int_\R v_1|v|^2(f-M_{\rm loc})\,\dd v.
\]
The extra \(|v|^2\) weights fast particles more heavily than
the momentum integral. Thus the odd part of \(f\) in \(v_1\)
can carry energy even when the momentum is zero.
For Figure~\ref{fig:stationary-velocity}, write \(\bar f\) for the average of
\(f\) over a small interval around \(x_1=1/4\). For more details, see
Appendix~\ref{app:numerical-illustrations}. The lower panels
show the current carried by particles with \(|v_1|\le V\):
\begin{equation}\label{eq:cumulative-energy-current}
 q_1(V)=\frac12\int_{|v_1|\le V}v_1|v|^2\bar f(v)\,\dd v.
\end{equation}

\subsection{A Gaussian Lower Bound}
\label{sec:stationary-bounds}

The Gaussian upper and lower bounds give
\(|\log F(x,v)|\le C(1+|v|^2)\),
so the entropy integrals are finite.

\begin{lemma}
\label{lem:stationary-gaussian-lower}
Let \(0\le\alpha<1\), and let
\(F\in L^\infty(W)\) be a nonnegative
stationary distributional solution. Then
\[
  F(x,v)\ge c e^{-C|v|^2}
  \qquad\text{for a.e.\ }(x,v)
\]
for some \(c,C>0\).
The constants may depend on \(\alpha\), \(\|F\|_{L^\infty(W)}\)
and the model parameters.
\end{lemma}
\begin{proof}
We first bound \(F\) from below on a velocity ball. Then
\eqref{eq:bath-kernel-lower} and Duhamel's formula give the
Gaussian lower bound at all velocities.

By \eqref{eq:collision-product-moment}, the gain and loss terms
belong to \(L^1(\T\times\R)\). The
\hyperref[par:mild-distributional]{characteristic argument}
therefore gives \eqref{eq:boltzmann-mild} for this stationary state.
The frequency bounds \eqref{eq:hard-frequency-comparison} and
\eqref{eq:collision-frequency-upper} give
\begin{equation}\label{eq:stationary-loss-bound}
  (1-\alpha)\nu[M_{\tau(x)}](v)+\alpha\nu[F](x,v)
  \le \Lambda\langle v\rangle^\gamma.
\end{equation}
Choose \(R\) so that
\(\|F\|_{L^\infty(W)}\int_{|v|>R}e^{-a|v|^2}\dd v\le 1/2\). Then
\[
  \int_{\T}\int_{|v|\le R}F(x,v)\,\dd v\dd x\ge\frac12.
\]

In \eqref{eq:boltzmann-mild}, retain just the paths with two bath
collisions and the bounded velocities used in
Proposition~\ref{prop:linear-minorisation}. Enlarge \(R\) to contain
these velocities. The loss rates are at most
\(\Lambda\langle R\rangle^\gamma\), and the gains carry the factor
\((1-\alpha)^2\). The same transport calculation, applied to the
mass of at least \(1/2\) above, gives constants \(c,r>0\) with
\[
 F(x,u)\ge c\one_{\{|u|\le r\}}
 \qquad\text{for a.e.\ }(x,u).
\]

Integrating \eqref{eq:bath-kernel-lower} over this velocity ball gives
\begin{equation}\label{eq:stationary-gain-lower}
 Q^+(F,M_{\tau(x)})(x,v)
 \geq c\langle v\rangle^{\gamma-2}M_{\tau(x)}(v).
\end{equation}%

Retaining this gain in \eqref{eq:boltzmann-mild} over a time
interval of length one and using \eqref{eq:stationary-loss-bound} gives
\[
\begin{aligned}
 F(x,v)&\geq c\langle v\rangle^{\gamma-2}e^{-|v|^2/(2\tau_-)}
       \int_0^1e^{-\Lambda\langle v\rangle^\gamma s}\,\dd s\\
 &\geq c\langle v\rangle^{-2}e^{-|v|^2/(2\tau_-)}
 \geq c'e^{-C|v|^2},
\end{aligned}
\]%
since the time integral is at least
\((1-e^{-\Lambda})/(\Lambda\langle v\rangle^\gamma)\).
\end{proof}

\subsection{Stationary Entropy Production}
\label{sec:stationary-entropy-3d}

The reservoir supplies energy to the gas at rate
\(\nabla_x\cdot q_f\), by \eqref{eq:stationary-energy-balance},
and its entropy changes at rate \(-\tau^{-1}\nabla_x\cdot q_f\).
At stationarity the gas entropy is constant. Thus the total
entropy production is \(-\int_\T\tau^{-1}\nabla_x\cdot q_f\,\dd x\),
which we relate to the collision dissipations:
\begin{align*}
 \mathcal D_{\rm bath}[f]
 &=-\int_{\T\times\R}Q(f,M_{\tau(x)})
                      \log\frac f{M_{\tau(x)}}\,\dd v\dd x,\\
 \mathcal D_{\rm self}[f]
 &=-\int_{\T\times\R}Q(f,f)\log f\,\dd v\dd x.
\end{align*}
For the homogeneous balance see \cite{CarlenLebowitzMouhot2015};
for the reservoir dissipation see
\cite[equations~(1.17)--(1.20)]{BisiCanizoLods2015}.

\begin{proposition}[Stationary entropy balance]
\label{prop:stationary-entropy-balance}
Let \(\gamma\geq0\), \(0\leq\alpha<1\), and let \(f\) be a
nonnegative stationary distributional solution of
\eqref{eq:boltzmann}. Assume that \(\tau\) satisfies
\eqref{eq:tau-bounds} and that, for some \(c,C>0\),
\[
 f(x,v)\leq Ce^{-c|v|^2}.
\]
Then \(\nabla_x\cdot q_f\in L^\infty(\T)\) and
\begin{equation}\label{eq:stationary-entropy-3d}
 -\int_\T\frac{\nabla_x\cdot q_f}{\tau}\,\dd x
 =(1-\alpha)\mathcal D_{\rm bath}[f]
       +\alpha\mathcal D_{\rm self}[f]\geq0.
\end{equation}%
Equality holds precisely when \(f=\rho M_\tau\) with constant
\(\rho\) and \(\tau\).  If
\(\tau^{-1}\in W^{1,\infty}(\T)\), then the left side equals
\(\int_\T q_f\cdot\nabla_x(\tau^{-1})\,\dd x\).
For a nonconstant reservoir, the entropy production is strictly positive,
and the energy current cannot have identically zero divergence.
\end{proposition}

\begin{proof}
Apply Lemma~\ref{lem:stationary-gaussian-lower} with
\(0<a<\min\{c,(2\tau_+)^{-1}\}\) to obtain a Gaussian lower bound for \(f\).
Testing the stationary equation with \(|v|^2/2\), which the
Gaussian bounds allow, gives
\begin{equation}\label{eq:stationary-energy-balance}
 \nabla_x\cdot q_f
 =\frac{1-\alpha}{2}\int_\R |v|^2Q(f,M_{\tau(x)})\,\dd v
 \in L^\infty(\T).
\end{equation}%
Write
\(Q[f]=(1-\alpha)Q(f,M_{\tau(x)})+\alpha Q(f,f)\).
The Gaussian bounds, \eqref{eq:strong-gain-form} and energy
conservation give \(|\log f|\leq C(1+|v|^2)\) and
\(|Q[f]|\leq C\langle v\rangle^\gamma e^{-c'|v|^2}\) for some \(c'>0\).
The stationary equation, chain rule and periodicity give
\[
 \int_{\T\times\R}Q[f]\log f\,\dd v\dd x
 =\int_{\T\times\R}v\cdot\nabla_x(f\log f-f)\,\dd v\dd x=0.
\]%
For distributional \(f\), take \(f_\varepsilon=\eta_\varepsilon*_xf\)
with a nonnegative periodic mollifier of mass one. Then
\(v\cdot\nabla_xf_\varepsilon=\eta_\varepsilon*_xQ[f]\).
Both Gaussian bounds persist under convolution, so Gaussian domination
gives \((\eta_\varepsilon*_xQ[f])\log f_\varepsilon\to Q[f]\log f\)
in \(L^1\).
Expanding \(Q[f]\), the self-collision contribution is
\(-\alpha\mathcal D_{\rm self}[f]\). For the reservoir term,
write \(\log f=\log(f/M_{\tau(x)})+\log M_{\tau(x)}\). Since
\(\log M_{\tau(x)}\) is a linear combination of \(1\) and \(|v|^2\),
\eqref{eq:stationary-energy-balance} gives
\[
 (1-\alpha)\int_\R Q(f,M_{\tau(x)})\log M_{\tau(x)}\,\dd v
 =-\frac{\nabla_x\cdot q_f}{\tau(x)},
\]%
which proves the identity in \eqref{eq:stationary-entropy-3d}. Put
\(M=M_{\tau(x)}\). Taking \(\psi=\log(f/M)\) and \(\psi=\log f\) in
\eqref{eq:symmetrised-weak-form}, and using
\eqref{eq:maxwellian-collision-balance} in the first case, gives
\begin{equation}\label{eq:dissipation-symmetrised}
\begin{aligned}
 \mathcal D_{\rm bath}[f]
 &=\frac12\int (fM_*-f'M_*')
             \log\frac{fM_*}{f'M_*'}\,B_n\,\dd n\dd v_*\dd v\dd x\geq0,\\
 \mathcal D_{\rm self}[f]
 &=\frac14\int (ff_*-f'f_*')
             \log\frac{ff_*}{f'f_*'}\,B_n\,\dd n\dd v_*\dd v\dd x\geq0.
\end{aligned}
\end{equation}%
If equality holds, the reservoir dissipation vanishes, and strict
positivity of the Carleman kernel on bounded velocity balls
(Lemma~\ref{lem:local-gain-lower}) forces \(f(x,v)=\rho(x)M_{\tau(x)}(v)\).
Both collision terms then vanish, so \(v\cdot\nabla_xf=0\). For a nonzero spatial Fourier frequency \(\xi\), the corresponding
coefficient of \(f\) is zero whenever \(\xi\cdot v\ne0\).
The exceptional plane has zero velocity measure, so this
coefficient vanishes almost everywhere. Hence \(f\), and with it
\(\rho\) and \(\tau\), are independent of \(x\). For constant \(\rho\) and \(\tau\), \(f=\rho M_\tau\) makes both
dissipations zero. The formula in terms of \(\nabla_x(\tau^{-1})\)
follows by periodic integration by parts.
\end{proof}

\subsection{Excluding Local Maxwellians}
\label{sec:entropy-general}
\label{sec:entropy-one-coordinate}

For a local Maxwellian, \(Q(F,F)=0\). We divide its stationary
equation by \(F\) and integrate in \(x\) to remove transport.
The gain at large velocities forces \(T\ge\tau\); the same
identity at \(v=0\) forces \(T=\tau\) and \(u=0\).

\begin{lemma}
\label{lem:stationary-maxwellian-rigidity}
Let \(\gamma\geq0\), \(0\leq\alpha<1\), and let \(\tau\) satisfy
\eqref{eq:tau-bounds}.  Let \(F\in L^\infty(W)\) be a
nonnegative stationary distributional solution of
\eqref{eq:boltzmann}.  If \(Q(F,F)=0\), then \(\tau\) is almost
everywhere constant and \(F\) is the corresponding spatially homogeneous,
centred Maxwellian.
\end{lemma}

\begin{proof}
The Gaussian upper bound and
Lemma~\ref{lem:stationary-gaussian-lower} give
\begin{equation}\label{eq:rigidity-two-sided}
 c e^{-C|v|^2}\leq F(x,v)\leq C e^{-c|v|^2}.
\end{equation}
By the equality case of the Boltzmann \(H\)-theorem
\cite[Sections~3.1--3.2 and~3.4]{CercignaniIllnerPulvirenti1994},
the assumption \(Q(F,F)=0\) implies
\[
 F(x,v)=\rho(x)M_{u(x),T(x)}(v)
 \quad\text{for almost every }x.
\]
The two Gaussian bounds imply
\begin{equation}\label{eq:rigidity-parameter-bounds}
 c\leq\rho\leq C,\qquad
 c\leq T\leq C,\qquad |u|\leq C.
\end{equation}%
Integrating the Gaussian bounds controls \(\rho\), \(|u|\) and \(T\)
from above and \(\rho\) from below. The upper bound at \(v=u\),
\(\rho(2\pi T)^{-3/2}\leq C\), then gives the lower bound on \(T\).

To remove transport after spatial integration, we divide the
stationary equation by \(F\). Since \(Q(F,F)=0\), the chain rule gives
\[
 (1-\alpha)\frac{Q(F,M_\tau)}{F}
 =\frac{v\cdot\nabla_xF}{F}
 =v\cdot\nabla_x\log F.
\]
For fixed \(v\), the spatial integral of the last term is zero
by periodicity. Since \(\alpha<1\),
\begin{equation}\label{eq:rigidity-averaged-collision}
 \int_\T\frac{Q(F,M_\tau)(x,v)}{F(x,v)}\,\dd x=0
 \qquad\text{for almost every }v\in\R.
\end{equation}
As in Proposition~\ref{prop:stationary-entropy-balance}, spatial
mollification justifies the chain rule. On bounded velocity sets,
\(F\) is bounded away from zero and \(Q(F,M_\tau)\) is bounded, by
\eqref{eq:rigidity-two-sided} and Lemma~\ref{lem:gain-bounds}.
The Maxwellian form of \(F\) and the bounds
\eqref{eq:rigidity-parameter-bounds} make the integral in
\eqref{eq:rigidity-averaged-collision} continuous in \(v\), by dominated
convergence in the gain formula \eqref{eq:strong-gain-form}.
Thus \eqref{eq:rigidity-averaged-collision} holds at every velocity,
and in particular,
\begin{equation}\label{eq:rigidity-gain-ratio}
 \int_\T\frac{Q^+(F,M_\tau)(x,v)}{F(x,v)}\,\dd x
 =\int_\T\nu[M_{\tau(x)}](v)\,\dd x
 \leq C\langle v\rangle^\gamma.
\end{equation}

The gain lower bound \eqref{eq:stationary-gain-lower} contains
\(M_{\tau(x)}\). If \(T<\tau\), this Gaussian decays more slowly than
\(F=\rho M_{u,T}\), contradicting the polynomial bound
\eqref{eq:rigidity-gain-ratio}. Dividing
\eqref{eq:stationary-gain-lower} by \(F\) and using
\eqref{eq:rigidity-parameter-bounds} gives, for \(R=|v|\geq1\),
\[
 \frac{Q^+(F,M_\tau)(x,v)}{F(x,v)}
 \geq cR^{\gamma-2}
 \exp\left[\left(\frac1{2T(x)}-\frac1{2\tau(x)}\right)R^2-CR\right].
\]
If \(\{T<\tau\}\) had positive measure, then
\(1/(2T)-1/(2\tau)\geq d\) on a set of positive measure, for some \(d>0\),
and integrating over this set would contradict
\eqref{eq:rigidity-gain-ratio} as \(R\to\infty\).  Hence
\[
 T\geq\tau\quad\text{almost everywhere}.
\]

Now evaluate the collision term at \(v=0\). Fix a spatial point where \(T\geq\tau\),
and write
\(s=v_*\cdot n\) and
\(\beta=(2\tau)^{-1}-(2T)^{-1}\geq0\).
At \(v=0\), the reflected velocities are
\(v'=sn\) and \(v_*'=v_*-sn\).  Direct substitution gives
\[
 \frac{F(sn)M_\tau(v_*-sn)}{F(0)M_\tau(v_*)}
 =\exp\left(\beta s^2+\frac{s\,u\cdot n}{T}\right).
\]%
By \eqref{eq:strong-gain-form} and
\(\nu[M_\tau](0)=\int B_nM_\tau(v_*)\,\dd n\dd v_*\), averaged under
\(v_*\mapsto-v_*\),
\[
 \frac{Q(F,M_\tau)(0)}{F(0)}
 =\int_{\R\times\mathbb S^2}B_n(-v_*,n)M_\tau(v_*)
 \left[e^{\beta s^2}
       \cosh\left(\frac{s\,u\cdot n}{T}\right)-1\right]
 \,\dd n\dd v_*.
\]%
Since \(\beta\geq0\) and \(\cosh\geq1\), the integrand is
nonnegative. It is strictly positive on a set of positive measure
whenever \(T>\tau\) or \(u\neq0\).
Equation~\eqref{eq:rigidity-averaged-collision} at \(v=0\) makes the
spatial integral zero, so
\[
 T=\tau,\qquad u=0\quad\text{almost everywhere}.
\]

Now \(F=\rho(x)M_{\tau(x)}\).  Both collision terms vanish,
so \(v\cdot\nabla_xF=0\). As in the proof of
Proposition~\ref{prop:stationary-entropy-balance}, this forces \(F\) to be
independent of \(x\), and its density and temperature are constant.
\end{proof}

\begin{remark}\label{rem:stationary-tails}
The same logarithmic chain rule applied to the stationary equation,
without assuming \(Q(F,F)=0\), gives
for almost every \(v\),
\[
 \int_\T\frac{(1-\alpha)Q(F,M_\tau)+\alpha Q(F,F)}{F}\,\dd x=0.
\]
Dropping the nonnegative self-collision gain and using
\eqref{eq:stationary-loss-bound} and~\eqref{eq:stationary-gain-lower},
\begin{equation}\label{eq:stationary-tail-bound}
 \int_\T\frac{M_{\tau(x)}(v)}{F(x,v)}\,\dd x\leq C\langle v\rangle^2
 \qquad\text{for almost every }v.
\end{equation}
In particular, a bound \(F\leq Ce^{-a|v|^2}\) forces
\(a\leq1/(2\operatorname{ess\,sup}\tau)\), so the restriction
\(a<1/(2\tau_+)\) on the weight \(W\) is sharp, up to the endpoint, when
\(\tau_+=\operatorname{ess\,sup}\tau\).

For Maxwell molecules, Proposition~\ref{prop:balance} and the Gaussian
lower bound in Lemma~\ref{lem:stationary-gaussian-lower} give
\eqref{eq:stationary-temperature-deficit} for \(F\). If \(\tau\) is not
almost everywhere constant, then \(T_F<\tau\) on a set of positive measure:
otherwise \(T_F=\tau\) and \(m_F=0\), so \(\nabla_x\cdot q_F=0\) by
\eqref{eq:balance-laws}, and Proposition~\ref{prop:stationary-entropy-balance}
would force \(\tau\) to be constant. The same Gaussian comparison in
\eqref{eq:stationary-tail-bound} rules out a bound
\(F\leq C\rho_FM_{u_F,T_F}\) on this set for any constant \(C\).
For Maxwell molecules this also shortens the
proof of the lemma: if \(F=\rho M_{u,T}\), then
\eqref{eq:stationary-tail-bound} gives \(T\geq\tau\), and
\eqref{eq:stationary-temperature-deficit} gives \(T=\tau\) and \(u=0\).
\end{remark}

\section{Relative Entropy Contraction and Initial Growth}
\label{sec:entropy}

When \(\alpha=0\), relative entropy to \(F_0\) is nonincreasing.
For small \(\alpha>0\) and a nonconstant reservoir, arbitrarily
small perturbations of \(F_\alpha\) can initially increase relative
entropy. The resulting solutions still converge exponentially by
Theorem~\ref{thm:boltzmann-contraction}.

For probability densities \(f\) and \(F>0\), the terms
\(-f+F\) in the introductory definition integrate to zero, giving
\[
 H(f\,|\,F)=\int_{\T}\int_{\R}f\log\frac fF\,\dd v\dd x,
\]
with \(0\log0=0\), allowing the value \(+\infty\).

\subsection{Contraction Under the Linear Flow \texorpdfstring{$(\alpha = 0)$}{(alpha = 0)}}
\label{sec:entropy-linear}

Applying the same Markov evolution to two probability densities
cannot increase their relative entropy. This is the data-processing
inequality \cite[Definition~2.1(iii), p.~303, and (3.1), p.~308]{Csiszar1967},
which gives, for probability densities \(f,g\),
\begin{equation}\label{eq:markov-relative-entropy-contraction}
 H(P_t^\tau f\,|\,P_t^\tau g)\leq H(f\,|\,g),
 \qquad t\geq0.
\end{equation}
We take \(g=F_0\), which the bath evolution leaves fixed.

\begin{proposition}[Contraction of relative entropy for the linear bath flow]
\label{prop:linear-relative-entropy}
Let \(\alpha=0\), and let \(F_0\) be the stationary state from
Theorem~\ref{thm:boltzmann-contraction}. Along every semigroup solution
\(f(t)=P_t^\tau f_{\rm in}\) with probability density \(f_{\rm in}\),
\[
  t\longmapsto H\bigl(f(t)\,\big|\,F_0\bigr)
\]
is nonincreasing, allowing the value \(+\infty\).
\end{proposition}

\begin{proof}
Apply \eqref{eq:markov-relative-entropy-contraction} at time \(t-s\)
to \(f(s)\) and \(F_0\). Since \(P_{t-s}^\tau F_0=F_0\),
\[
 H(f(t)\,|\,F_0)\leq H(f(s)\,|\,F_0),\qquad 0\leq s\leq t.
\]
\end{proof}

\subsection{Initial Growth of Relative Entropy}
\label{sec:entropy-growth}

For a nonconstant reservoir, \(F\) is not a local Maxwellian by
Lemma~\ref{lem:stationary-maxwellian-rigidity}. Hence
\(F'F_*'-FF_*\) does not vanish identically. We use this to choose initial data whose solution has relative entropy
greater than its initial value for small positive times.

Let \(F>0\) be a stationary state in
\(L^1(W)\cap L^\infty(W)\).  Whenever the
transport derivative is integrable and the pairing is finite, define
\begin{equation}\label{eq:relative-entropy-pairing}
 \mathcal E_F(f)
 =\int_{\T\times\R}\log\frac fF
 \bigl[-v\cdot\nabla_xf+(1-\alpha)Q(f,M_{\tau(x)})+\alpha Q(f,f)\bigr]
 \,\dd v\dd x.
\end{equation}%
For smooth positive solutions,
\(\mathcal E_F(f)=\frac{\dd}{\dd t}H(f\,|\,F)\). For the mild solutions used
below, \eqref{eq:entropy-initial-lower-bound} gives
\(\liminf_{t\downarrow0}t^{-1}[H(f(t)\,|\,F)-H(f_{\rm in}\,|\,F)]
\geq\mathcal E_F(f_{\rm in})\).

All quantities at a collision are evaluated at the same spatial point.
Write \(r=f/F\), \(M=M_{\tau(x)}\) and \(\Delta_F=F'F_*'-FF_*\).
Collision integrals written without a measure are taken with respect to
\(B_n\,\dd n\dd v_*\dd v\dd x\), which is invariant under the changes of
variables used in \eqref{eq:symmetrised-weak-form}.

\begin{proposition}[Rate of change of relative entropy]
\label{prop:entropy-three-terms}
Suppose that \(0<c\leq r\leq C<\infty\) and
\(v\cdot\nabla_xf\in L^1\).  Then
\begin{equation}\label{eq:entropy-three-terms}
\begin{aligned}
 \mathcal E_F(f)
 ={}&-(1-\alpha)\int FM_*\Bigl[r\log\frac r{r'}-r+r'\Bigr]
  -\frac\alpha2\int FF_*\Bigl[rr_*\log\frac{rr_*}{r'r_*'}-rr_*+r'r_*'\Bigr]\\
  &+\frac\alpha2\int\Delta_F(r-1)(r_*-1).
\end{aligned}
\end{equation}%
The first two terms are nonpositive by
\(s\log(s/t)-s+t\geq0\). The last term can have either sign.
\end{proposition}

\begin{proof}
Apply the chain rule and periodicity to \(F(r\log r-r+1)\):
\[
 -\int\log r\,v\cdot\nabla_xf\,\dd v\dd x
 =\int(1-r)\,v\cdot\nabla_xF\,\dd v\dd x
 =-\int r\,v\cdot\nabla_xF\,\dd v\dd x.
\]%
As in Proposition~\ref{prop:stationary-entropy-balance}, justify this
by mollifying \(f\) and \(F\) in \(x\) with the same kernel. This
preserves \(c\leq f_\varepsilon/F_\varepsilon\leq C\); since
\(v\cdot\nabla_xf,v\cdot\nabla_xF\in L^1\), the mollified transport
derivatives converge in \(L^1\), and the bounded ratio factors allow
passage to the limit.
Substituting the stationary equation for \(F\), the
reservoir contribution becomes
\[
 (1-\alpha)\int FM_*\Bigl[r\log\frac{r'}r-(r'-r)\Bigr],
\]%
the first term of \eqref{eq:entropy-three-terms}.
Exchanging \(v\) and \(v_*\), as in \eqref{eq:symmetrised-weak-form},
writes the self-collision contribution as
\[
 \frac\alpha2\int FF_*
 \left[rr_*\log\frac{r'r_*'}{rr_*}
       -(r'+r_*'-r-r_*)\right].
\]%
Adding and subtracting \(r'r_*'-rr_*\) separates this into
the second term of \eqref{eq:entropy-three-terms} and
\[
 \frac\alpha2\int FF_*
 \left[(r'-1)(r_*'-1)-(r-1)(r_*-1)\right].
\]%
The pre/post-collisional change of variables in the first summand gives
the third term.
\end{proof}

If \(F\) is a local Maxwellian, then \(\Delta_F=0\) and
\(\mathcal E_F\leq0\). The two dissipations in \eqref{eq:entropy-three-terms}
are taken relative to \(F\), unlike those of
Proposition~\ref{prop:stationary-entropy-balance}.

We now expand \(\mathcal E_F\) around \(F\). Let
\begin{equation}\label{eq:entropy-linearised-operator}
 \mathcal L_Fg=-v\cdot\nabla_xg+(1-\alpha)Q(g,M_\tau)
                 +\alpha\{Q(g,F)+Q(F,g)\}
\end{equation}%
be the linearisation of the equation at \(F\), let \(h\) be a
bounded real function with \(\int Fh\,\dd v\dd x=0\) and
\(v\cdot\nabla_x(Fh)\in L^1\), and set \(f_\varepsilon=F(1+\varepsilon h)\)
for \(|\varepsilon|\|h\|_\infty\leq1/2\). 
Stationarity and bilinearity make the expression in brackets in
\eqref{eq:relative-entropy-pairing} equal to
\(\varepsilon\mathcal L_F(Fh)+\alpha\varepsilon^2Q(Fh,Fh)\).
Both terms are integrable by the transport assumption, the Gaussian
bounds and the polynomial collision kernel. The uniform expansion
\(\log(1+\varepsilon h)=\varepsilon h+O(\varepsilon^2)\) therefore gives
\begin{equation}\label{eq:general-production-expansion}
 \mathcal E_F(f_\varepsilon)
 =\varepsilon^2\int h\mathcal L_F(Fh)\,\dd v\dd x+O(\varepsilon^3).
\end{equation}%

Expanding \eqref{eq:entropy-three-terms} to second order in
\(\varepsilon\), justified by the same bounds, and comparing with
\eqref{eq:general-production-expansion} gives
\begin{equation}\label{eq:linearised-three-terms}
\begin{aligned}
 \int h\mathcal L_F(Fh)\,\dd v\dd x
 ={}&-\frac{1-\alpha}2\int FM_*(h'-h)^2
  -\frac\alpha4\int FF_*(h'+h_*'-h-h_*)^2\\
  &+\frac\alpha2\int\Delta_Fhh_*.
\end{aligned}
\end{equation}%

A density perturbation, whose ratio to \(F\) depends only on \(x\),
cancels in both collision differences in
\eqref{eq:linearised-three-terms}. We add a small perturbation of
the velocity distribution. The two dissipations are quadratic in
its size, whereas its interaction with the density perturbation
has a linear term. Choosing its sign and size makes the total
positive. We can then scale the whole perturbation down to make
the initial datum as close to \(F\) as desired.

\begin{theorem}[Initial entropy growth]
\label{thm:ness-entropy-growth}
Fix \(A,S<\infty\) and \(H>0\). Take \(0<\alpha\le1/2\) small
enough for Lemma~\ref{thm:boltzmann-local} on \([0,H]\) with these
bounds, and let
\(F\in L^1(W)\cap L^\infty(W)\) be a nonnegative
stationary distributional solution, with
\[
 \norm{F}_{L^1(W)}<A,
 \qquad \norm{F}_{L^\infty(W)}<S,
\]
and with
\(Q(F,F)\not\equiv0\). For every \(\delta>0\), there is a probability density
\(f_{\rm in}\in L^1(W)\cap L^\infty(W)\) such that
\[
 \norm{f_{\rm in}-F}_{L^1(W)}+\norm{f_{\rm in}-F}_{L^\infty(W)}<\delta,
 \qquad
 \frac12\le \frac{f_{\rm in}}{F}\le\frac32\quad\text{a.e.}
\]
We may also require
\(\norm{f_{\rm in}}_{L^1(W)}\le A\) and
\(\norm{f_{\rm in}}_{L^\infty(W)}\le S\). The mild solution
on \([0,H]\) then satisfies, for some \(0<t_0\le H\),
\[
  H(f(t)\,|\,F)>H(f_{\rm in}\,|\,F),
  \qquad 0<t<t_0.
\]
\end{theorem}

\begin{proof}[Proof of Theorem~\ref{thm:ness-entropy-growth}]
The ratio \(Q(F,F)/F\) is locally integrable by the Gaussian bounds
and collision estimates. Choose a proper spatial ball \(U\Subset\T\) on
which \(Q(F,F)\not\equiv0\), and fix \(\chi\in C_c^\infty(\R)\)
with \(\int\chi\,\dd v=1\).  For
\(\psi\in C_c^\infty(U\times\R)\), put
\[
 \zeta(x,v)=\psi(x,v)-\chi(v)\int_\R\psi(x,w)\,\dd w.
\]
Then \(\zeta\) is smooth, compactly supported and of zero local mass.
If \(\int\zeta\,Q(F,F)/F=0\) for every such \(\psi\), then
\(Q(F,F)/F\) would be independent of \(v\) on \(U\times\R\); collision mass
conservation, \(\int Q(F,F)\,\dd v=0\), and \(\rho_F>0\) would then force it
to vanish there, a contradiction.  Choose \(\zeta\), changing its sign if
necessary, so that
\[
 J=\int_{\T\times\R}\frac\zeta FQ(F,F)\,\dd v\dd x>0.
\]%
Since the spatial support is a proper compact subset of the torus and
\(\rho_F>0\), there is \(\eta\in C^\infty(\T)\), equal to one on
that support, such that \(\int_\T\eta\rho_F\,\dd x=0\).
By Lemma~\ref{lem:stationary-gaussian-lower}, \(\zeta/F\) is
bounded. For \(k>0\), set \(g=\eta F+k\zeta\).
The term \(\eta F\) changes only the local density, while
\(\zeta\) has zero mass at every \(x\). Both have zero total mass.
Moreover,
\[
 v\cdot\nabla_xg
 =Fv\cdot\nabla_x\eta+\eta v\cdot\nabla_xF
   +kv\cdot\nabla_x\zeta\in L^1,
\]
since \(v\cdot\nabla_xF\) is the stationary collision term.
Thus \(h=g/F\) satisfies the hypotheses of
\eqref{eq:linearised-three-terms}. Appendix~\ref{app:entropy-growth-proof}
will justify the initial entropy bound for the mild solution.

In \eqref{eq:linearised-three-terms}, the part \(\eta\) of
\(h=\eta+k\zeta/F\) cancels from both collision differences.
Both dissipations therefore contain \(k^2\). For the last term,
collision mass conservation and \(\eta=1\) on the support of
\(\zeta\) give
\[
 \int\Delta_F\eta^2=0,\qquad
 \int\Delta_F\eta\,\frac{\zeta}{F}
 =\int\eta\,\frac{\zeta}{F}Q(F,F)\,\dd v\dd x=J.
\]%
Exchanging \(v\) and \(v_*\) gives the same cross term a second
time. Hence
\[
 \int\frac gF\mathcal L_Fg\,\dd v\dd x=\alpha kJ+k^2D,\qquad
 D=\int\frac{\zeta}{F}\mathcal L_F\zeta\,\dd v\dd x.
\]%
Since \(J>0\), we may choose \(k>0\) with \(\alpha kJ+k^2D>0\); fix
this \(k\). Then \eqref{eq:general-production-expansion} gives
\[
 \mathcal E_F(F+\varepsilon g)
 =\varepsilon^2(\alpha kJ+k^2D)+O(\varepsilon^3)>0
\]
for all sufficiently small \(\varepsilon>0\).

The datum \(f_{\rm in}=F+\varepsilon g\) satisfies
\(\int_{\T\times\R}f_{\rm in}=1\) and
\[
 \|f_{\rm in}-F\|_{L^1(W)}+
 \|f_{\rm in}-F\|_{L^\infty(W)}\leq C\varepsilon.
\]
Take \(\varepsilon\) so small that this is less than \(\delta\).
Since \(\|F\|_{L^1(W)}<A\),
\(\|F\|_{L^\infty(W)}<S\), and \(g/F\) is bounded, decrease
\(\varepsilon\) further to obtain
\[
 \|f_{\rm in}\|_{L^1(W)}<A,\qquad
 \|f_{\rm in}\|_{L^\infty(W)}<S,\qquad
 \frac12\leq\frac{f_{\rm in}}F\leq\frac32.
\]
Lemma~\ref{thm:boltzmann-local} gives the local solution.
Its strict initial entropy growth follows from
\eqref{eq:entropy-initial-lower-bound} and
\(\mathcal E_F(f_{\rm in})>0\).
\end{proof}

\subsection{Entropy Growth and Return to the Stationary State}
\label{sec:entropy-selected}
\label{sec:entropy-return}

\begin{corollary}
\label{cor:ness-entropy-coexistence}
Assume \(0<\alpha<\alpha_0\) in
Theorem~\ref{thm:boltzmann-contraction}, let \(F=F_\alpha\), and suppose
that \(\tau\in L^\infty(\T)\) is not almost everywhere constant.
Then \(Q(F,F)\not\equiv0\).  For every \(\delta>0\), the datum in
Theorem~\ref{thm:ness-entropy-growth} may be chosen in \(\mathcal R\),
with
\[
 \|f_{\rm in}-F\|_{L^1(W)}+\|f_{\rm in}-F\|_{L^\infty(W)}<\delta.
\]
Its solution is global and satisfies the exponential return estimate
in Theorem~\ref{thm:intro-relaxation}. Its relative entropy to \(F\)
exceeds its initial value for all sufficiently small positive times.
\end{corollary}

\begin{proof}
Lemma~\ref{lem:stationary-maxwellian-rigidity} gives \(Q(F,F)\not\equiv0\).
Applying \eqref{eq:interval-return} to the stationary solution gives
\[
 \|F\|_{L^1(W)}+\|F\|_{L^\infty(W)}\leq\frac{3R}{4}.
\]
Apply Theorem~\ref{thm:ness-entropy-growth} with
\(A=S=R\) and \(H=t_c\), choosing the sum of the two distances
smaller than \(\min\{\delta,R/4\}\).
Then \(f_{\rm in}\in\mathcal R\). The solution from
Theorem~\ref{thm:boltzmann-contraction} agrees with the local
solution by Lemma~\ref{thm:boltzmann-local}, and therefore has
the initial entropy increase.
\end{proof}

Relative entropy still tends to zero exponentially. The Gaussian
bounds give an upper bound by a constant times the weighted \(L^1\)
distance to \(F_\alpha\).

\begin{corollary}[Exponential upper bound for relative entropy]
\label{cor:stationary-entropy-decay}
Under the hypotheses of Theorem~\ref{thm:boltzmann-contraction}, every solution
starting from \(f_{\rm in}\in\mathcal R\) satisfies
\[
\begin{aligned}
  \frac12\norm{f(t)-F_\alpha}_{L^1}^2
  &\le H\bigl(f(t)\,\big|\,F_\alpha\bigr)\\
  &\le C\,2^{-\lfloor t/t_c\rfloor}
       \bigl(\|f_{\rm in}-F_\alpha\|_{L^1(W)}
       +\|f_{\rm in}-F_\alpha\|_{L^\infty(W)}\bigr).
\end{aligned}
\]
\end{corollary}

\begin{proof}
Theorem~\ref{thm:boltzmann-contraction} and
Lemma~\ref{lem:stationary-gaussian-lower} give
\(f(t),F_\alpha\le Ce^{-a|v|^2}\) and \(F_\alpha\ge ce^{-C|v|^2}\),
uniformly for \(t\ge0\). For
\(\varphi(r)=r\log r-r+1\), the probability normalisation gives
\[
  H(f\,|\,F_\alpha)
  =\int_{\T}\int_{\R}F_\alpha\varphi(f/F_\alpha)\,\dd v\dd x.
\]
On \(\{f<F_\alpha\}\), \(F_\alpha\varphi(f/F_\alpha)\le F_\alpha-f\).  On
\(\{f\ge F_\alpha\}\),
\[
  F_\alpha\varphi(f/F_\alpha)
  \le(f-F_\alpha)\log(f/F_\alpha)
  \le C(1+|v|^2)|f-F_\alpha|.
\]
Since \(1+|v|^2\le C W(v)\), the entropy is bounded by
\(C\norm{f-F_\alpha}_{L^1(W)}\).
Theorem~\ref{thm:boltzmann-contraction} then gives the upper bound.
Pinsker's inequality gives the lower bound.

\end{proof}

\section*{Use of large language models}
We used large language models (LLMs) extensively in reviewing the
literature. We also used
them to rework proofs, look for mistakes throughout the manuscript, and
help with calculations, simulation code, typesetting, and editing.

Our initial plan was to look for stationary states by a fixed-point
argument. The mild-solution approach for \(\gamma=0\) grew out of this
work. Discussions with LLMs were particularly helpful in putting together
the hard-potential estimates and working through the details needed to
justify formal arguments in the entropy section.  We have
checked and revised the arguments and take responsibility for the paper.

\section*{Acknowledgments}
The authors would like to thank Angeliki Menegaki for many helpful discussions during the preparation of this paper. 
Both authors are supported by a Royal Society University Research Fellowship R1 251808.

 For the purpose of
open access, the author has applied a Creative Commons Attribution (CC-BY) licence to any Author
Accepted Manuscript version arising from this submission.
\appendix

\section{Initial Entropy Growth for Mild Solutions}
\label{app:entropy-growth-proof}
\label{app:relative-entropy-derivative}

Section~\ref{sec:entropy-growth} constructs a perturbation
whose relative entropy initially increases. Here we justify the initial
\(L^1\) derivative of its mild solution and use convexity to obtain
\eqref{eq:entropy-initial-lower-bound}.

Write the collision part of the equation as
\[
 \mathcal C(g)=(1-\alpha)Q(g,M_{\tau(x)})+\alpha Q(g,g).
\]
The gain and loss estimates give, on every set bounded in \(L^\infty(W)\),
\begin{equation}\label{eq:collision-moment-lipschitz}
 \norm{\mathcal C(f)-\mathcal C(g)}_{L^1}
 \le C\norm{f-g}_{L^1(\langle v\rangle^\gamma)}
 \le C'\norm{f-g}_{L^1(W)}.
\end{equation}
Apply \eqref{eq:collision-product-moment} to
\(Q(f-g,f)+Q(g,f-g)\) and \(Q(f-g,M_{\tau(x)})\), then integrate
in \(x\). The Gaussian bounds make the velocity moments of
\(f,g,M_\tau\) uniformly bounded in \(x\).
Finally, \(\langle v\rangle^\gamma\leq CW(v)\) gives the second inequality.

For a perturbation \(g=\eta F+z\), with \(\eta\in C^\infty(\T)\) and
\(z\in C_c^\infty(\T\times\R)\), the linear term in the equation at \(F\) is
\(\mathcal L_Fg\), with \(\mathcal L_F\) defined in
\eqref{eq:entropy-linearised-operator}.
The stationary equation gives
\[
 v\cdot\nabla_xg=Fv\cdot\nabla_x\eta+\eta\mathcal C(F)
                         +v\cdot\nabla_xz\in L^1.
\]
The Gaussian bounds and the collision estimate used in
\eqref{eq:collision-moment-lipschitz} also give
\(\mathcal L_Fg,Q(g,g)\in L^1\).
Assume \(\int g=0\), and put \(f_{\rm in}=F+\varepsilon g\).
Since \(z\) has bounded velocity support,
Lemma~\ref{lem:stationary-gaussian-lower} makes
\(g/F=\eta+z/F\) bounded. Thus, for small
\(\varepsilon>0\) the datum is a positive probability density and its ratio to
\(F\) is bounded above and below. The right side of its evolution equation is
\begin{equation}\label{eq:constructed-initial-derivative}
 R_0=\varepsilon\mathcal L_Fg+\alpha\varepsilon^2Q(g,g)\in L^1.
\end{equation}
In particular, \(-v\cdot\nabla_xf_{\rm in}\in L^1\).

Let \(T_t\) be the free transport semigroup of
Section~\ref{sec:linear-harris}. Since \(v\cdot\nabla_xf_{\rm in}\in L^1\),
\((T_tf_{\rm in}-f_{\rm in})/t\to-v\cdot\nabla_xf_{\rm in}\)
in \(L^1\) as \(t\downarrow0\). By the
\hyperref[par:mild-distributional]{method of characteristics}, for almost
every \((y,v)\) the function \(t\mapsto f(t,y+tv,v)\) is absolutely
continuous with derivative \(\mathcal C(f(t))(y+tv,v)\); integrating gives
\begin{equation}\label{eq:appendix-transport-duhamel}
 f(t)=T_tf_{\rm in}+\int_0^tT_{t-s}\mathcal C(f(s))\dd s.
\end{equation}
Subtract \(f_{\rm in}\) from \eqref{eq:appendix-transport-duhamel} and divide by
\(t\). Continuity of \(f\) in \(L^1(W)\),
\eqref{eq:collision-moment-lipschitz}, and the limit
\(\|T_t h-h\|_{L^1}\to0\) for \(h\in L^1\) give
\[
 \frac{f(t)-f_{\rm in}}{t}\longrightarrow
 -v\cdot\nabla_xf_{\rm in}+\mathcal C(f_{\rm in})=R_0
 \qquad\text{in }L^1.
\]

The Gaussian bounds make both entropies finite. Convexity of
\(r\log r-r+1\) and the probability normalisation give
\[
 H(f(t)\,|\,F)-H(f_{\rm in}\,|\,F)
 \geq\int_{\T\times\R}\log\frac{f_{\rm in}}F
                  (f(t)-f_{\rm in})\,\dd v\dd x.
\]
Since \(\log(f_{\rm in}/F)\) is bounded, the preceding \(L^1\) limit yields
\begin{equation}\label{eq:entropy-initial-lower-bound}
  \liminf_{t\downarrow0}
  \frac{H(f(t)\,|\,F)-H(f_{\rm in}\,|\,F)}{t}
  \geq\mathcal E_F(f_{\rm in}),
\end{equation}
with \(\mathcal E_F\) defined in \eqref{eq:relative-entropy-pairing}.

\section{Numerical Illustrations}
\label{app:numerical-illustrations}

The figures use \(\tau(x_1)=1+0.3\cos(2\pi x_1)\), unit mass and
\(\alpha=0.2\). The dynamical stability
of these states at \(\alpha=0.2\) remains to be determined.

\paragraph{Velocity-space plots.}
For Figure~\ref{fig:stationary-velocity}, \(\bar f(v)\) is the average of
\(f\) over the interval of width \(1/128\) centred at \(x_1=1/4\).
The top panels show the relative deviation of \(\bar f\)
from the Maxwellian with the same density, mean velocity and temperature
as \(\bar f\).

\paragraph{Models and collision rates.}
The reference Boltzmann kernels are \(B_\sigma=|v-v_*|^\gamma/(4\pi)\),
with uniform outgoing direction and \(\gamma=0,1\).
The centred BGK comparison solves
\[
 v_1\partial_{x_1}f
 =\nu\bigl[\alpha\rho_f M_{T_f}
          +(1-\alpha)\rho_f M_\tau-f\bigr].
\]
Comparing the models requires a choice of collision rates, and there is no
canonical way to match them. We choose the rates so that a small energy
current near the unit Maxwellian has the same initial decay rate in each
model.

Fix the bath at the Maxwellian \(M_1\) of density and temperature one.
A temperature gradient at constant pressure gives the streaming term
\[
 v_1\left.\partial_T(T^{-1}M_T)\right|_{T=1}
 =\tfrac12M_1v_1(|v|^2-5)=\tfrac12h.
\]
We linearise in this direction. The mass and energy moments of
\(h\) vanish because \(h\) is odd, and Gaussian integration gives
\(\int vh\,\dd v=0\) and \((q_h)_1=\frac12\int |v|^2v_1h\,\dd v=5\).

Differentiating the collision term gives
\(L_\gamma h=Q(h,M_1)+\alpha Q(M_1,h)\).
The initial decay rate is
\(\lambda_\gamma=-(q_{L_\gamma h})_1/(q_h)_1\).
At \(\alpha=1/5\), Gaussian integration of the bath and gas terms gives
\begin{equation}
 \lambda_0=\frac{4(10/3)+5/3}{25}=\frac35,
 \qquad
 \lambda_1=\frac{4(23)+32/3}{25\sqrt\pi}=\frac{308}{75\sqrt\pi}.
 \label{eq:numerical-current-rate}
\end{equation}
In the figures, we use \(\nu=3/5\) for BGK and multiply the hard sphere
kernel by \((3/5)/\lambda_1=45\sqrt\pi/308\).
By \eqref{eq:numerical-current-rate}, this matches the initial current decay
rate while retaining its dependence on relative velocity.

Details of the numerical method will be presented in separate
work.

{\fontsize{9}{11}\selectfont

}

\end{document}